\documentclass{article}
\usepackage{verbatim}
\usepackage{amsmath}
\usepackage[mathscr]{eucal}
\usepackage{amsthm}
\usepackage{amssymb,mathrsfs}
\usepackage[english,german]{babel}
\usepackage{graphicx}
\usepackage{array}
\theoremstyle{plain}
\newtheorem{Theo}{Theorem}[section]
\newtheorem{Prop}[Theo]{Proposition}

\newtheorem{Cor}[Theo]{Corollary}
\newtheorem{Main}{Theorem}
\theoremstyle{definition}
\newtheorem{Def}[Theo]{Definition}
\theoremstyle{remark}
\newtheorem{Rem}[Theo]{Remark}
\begin{document}

\selectlanguage{english}
\title{Involutions of Iwahori-Hecke algebras and representations of 
fixed subalgebras}
\author{Hideo Mitsuhashi \\
Department of Mechanical Engineering \\
Tomakomai National College of Technology \\
443 Nishikioka, Tomakomai, Hokkaido, 059-1275, Japan}

\date{}
\maketitle

\renewcommand{\theenumi}{\arabic{enumi}}
\renewcommand{\labelenumi}{(\theenumi)}

\begin{abstract}
We establish branching rules between some Iwahori-Hecke algebra of type $B$ and 
their subalgebras which are defined as fixed subalgebras by involutions 
including Goldman involution. 
The Iwahori-Hecke algebra of type $D$ is one of such fixed subalgebras. 
We also obtain branching rules between those fixed subalgebras 
and their intersection subalgebra. 
We determine basic sets of irreducible representations of those fixed 
subalgebras and their intersection subalgebra by making use of generalized Clifford theory. 
\end{abstract}

\section{Introduction}
In our paper \cite{Mitsuhashi}, we gave the $q$-analogue of the alternating group 
and its irreducible representations. 
Since the Iwahori-Hecke algebra of type $A$ can be considered as a $q$-analogue of the 
symmetric group, we defined the $q$-analogue of the alternating group as a subalgebra 
of it. 
An extension of \cite{Mitsuhashi} to the Hecke algebra of type $B$ was given in 
\cite{Mitsuhashi2}. 
We showed that those subalgebras can be defined as fixed subalgebras by 
Goldman involutions of Hecke algebras in \cite{Mitsuhashi2}. 
The Iwahori-Hecke algebra of type $D$ can be considered as a subalgebra of some 
specialized Iwahori-Hecke algebra of type $B$, and can be realized as 
a fixed subalgebra by an involution, which we denote by $\sharp$, 
of the specialized Iwahori-Hecke algebra. 
The specialized Iwahori-Hecke algebra has a $\mathbb{Z}_2$-graded 
Clifford system whose submodule including identity element coincides with 
the Iwahori-Hecke algebra of type $D$. 
Using generalized Clifford theory, we give a new proof of the branching rule between 
the specialized Iwahori-Hecke algebra 
and the Iwahori-Hecke algebra of type $D$, which has already been obtained in 
\cite{Hoefsmit}. 
We do not specify the base fields in the statements below in this section 
because of simplicity. The details are specified in corresponding paragraphs 
since the next chapter. 
Let $\pi_{\lambda}$ be the irreducible representation of the Iwahori-Hecke algebra of 
type $B$ corresponding to $\lambda=(\lambda^{(1)},\lambda^{(2)})$, a $2$-tuple of 
Young diagrams. 
Then we have the following.
\renewcommand{\theMain}{6.8}
\begin{Main}
\begin{enumerate}
\renewcommand{\theenumi}{\arabic{enumi}}
\renewcommand{\labelenumi}{($\sharp$\theenumi)}
\rm \item \it
The restrictions $\pi_{\lambda,\sharp}$ and $\pi_{\lambda^*,\sharp}$ of irreducible 
representations $\pi_{\lambda}$ and $\pi_{\lambda^*}$ corresponding to 
$\lambda=(\lambda^{(1)},\lambda^{(2)})$ and $\lambda^*=(\lambda^{(2)},\lambda^{(1)})$ 
to the fixed subalgebra by $\sharp$ are equivalent. 
\rm \item \it
If $\lambda^*=\lambda$, then 
$\pi_{\lambda,\sharp}$ decomposes into two inequivalent 
subrepresentations $\pi_{\lambda,\sharp}^+$ and $\pi_{\lambda,\sharp}^-$ of the same degree. 
\rm \item \it
Let $\{\lambda_i,\lambda^*_i,\mu_j\}_{i,j}$ be the set of $2$-tuples of Young diagrams such that 
$\lambda_i{\neq}\lambda^*_i$ and $\mu_j=\mu^*_j$. Then 
$\{\pi_{\lambda_i,\sharp},\pi_{\mu_j,\sharp}^+,\pi_{\mu_j,\sharp}^-\}_{i,j}$ is a basic set of 
irreducible representations of the fixed subalgebra by $\sharp$. 
\end{enumerate}
\end{Main}
Besides such involution, there are two other involutions of such specialized 
Iwahori-Hecke algebra. One of them is Goldman involution, which we denote by $\natural$. 
Goldman involution itself can be defined for generic (non-specialized) Iwahori-Hecke algebras 
of type $B$. 
We have already given the fixed subalgebra by Goldman involution in generic 
Iwahori-Hecke algebras of type $B$, and have determined the branching rule and 
a basic set of irreducible representations in \cite{Mitsuhashi2} as follows. 
\renewcommand{\theMain}{6.2}
\begin{Main}
\begin{enumerate}
\renewcommand{\theenumi}{\arabic{enumi}}
\renewcommand{\labelenumi}{($\natural$\theenumi)}
\rm \item \it
The restrictions $\pi_{\lambda,\natural}$ and $\pi_{\lambda'^*,\natural}$ of irreducible 
representations $\pi_{\lambda}$ and $\pi_{\lambda'^*}$ corresponding to 
$\lambda=(\lambda^{(1)},\lambda^{(2)})$ and $\lambda'^*=(\lambda^{(2)}{'},\lambda^{(1)}{'})$ 
to the fixed subalgebra by $\natural$ are equivalent, 
where $\lambda^{(k)}{'}$ stands for the transpose of $\lambda^{(k)}$. 
\rm \item \it
If $\lambda'^*=\lambda$, then 
$\pi_{\lambda,\natural}$ decomposes into two inequivalent 
subrepresentations $\pi_{\lambda,\natural}^+$ and $\pi_{\lambda,\natural}^-$ of the same degree. 
\rm \item \it
Let $\{\lambda_i,\lambda_i'^*,\mu_j\}_{i,j}$ be the set of $2$-tuples of Young diagrams such that 
$\lambda_i{\neq}\lambda_i'^*$ and $\mu_j=\mu_j'^*$. Then 
$\{\pi_{\lambda_i,\natural},\pi_{\mu_j,\natural}^+,\pi_{\mu_j,\natural}^-\}_{i,j}$ is a basic set of 
irreducible representations of the fixed subalgebra by $\natural$. 
\end{enumerate}
\end{Main}
The other, which we denote by $\flat$, is new. We give the fixed subalgebra by it 
and determine properties of it. 
In the same manner as in the case $\sharp$ or $\natural$, we obtain the branching rule 
and a basic set of irreducible representations of the fixed subalgebra. 
\renewcommand{\theMain}{6.13}
\begin{Main}
\begin{enumerate}
\renewcommand{\theenumi}{\arabic{enumi}}
\renewcommand{\labelenumi}{($\flat$\theenumi)}
\rm \item \it
The restrictions $\pi_{\lambda,\flat}$ and $\pi_{\lambda',\flat}$ of irreducible 
representations $\pi_{\lambda}$ and $\pi_{\lambda'}$ corresponding to 
$\lambda=(\lambda^{(1)},\lambda^{(2)})$ and $\lambda'=(\lambda^{(1)}{'},\lambda^{(2)}{'})$ 
to the fixed subalgebra by $\flat$ are equivalent. 
\rm \item \it
If $\lambda'=\lambda$, then 
$\pi_{\lambda,\flat}$ decomposes into two inequivalent 
subrepresentations $\pi_{\lambda,\flat}^+$ and $\pi_{\lambda,\flat}^-$ of the same degree. 
\rm \item \it
Let $\{\lambda_i,\lambda'_i,\mu_j\}_{i,j}$ be the set of $2$-tuples of Young diagrams such that 
$\lambda_i{\neq}\lambda'_i$ and $\mu_j=\mu'_j$. Then 
$\{\pi_{\lambda_i,\flat},\pi_{\mu_j,\flat}^+,\pi_{\mu_j,\flat}^-\}_{i,j}$ is a basic set of 
irreducible representations of the fixed subalgebra by $\flat$. 
\end{enumerate}
\end{Main}
The intersection of those three fixed subalgebras can be considered as the 
fixed subalgebra by Goldman involution in the Iwahori-Hecke algebra of type $D$, 
while it may be viewed as a subalgebra of both the fixed subalgebra by $\natural$ 
and that by $\flat$. 
We denote by $\pi_{\lambda,\dagger}$, 
$\pi_{\lambda,\natural,\dagger}$, 
$\pi_{\lambda,\sharp,\dagger}$, 
$\pi_{\lambda,\flat,\dagger}$ the restrictions of 
$\pi_\lambda$,$\pi_{\lambda,\natural}$,
$\pi_{\lambda,\sharp}$,
$\pi_{\lambda,\flat}$ to the intersection respectively. 
It is clear that 
$\pi_{\lambda,\dagger}{\cong}\pi_{\lambda,\natural,\dagger}{\cong}
\pi_{\lambda,\sharp,\dagger}{\cong}\pi_{\lambda,\flat,\dagger}$. 
We give branching rules between those fixed subalgebras and the intersection subalgebra of them, 
and determine a basic set of irreducible representations of the intersection subalgebra as follows. 
\begin{enumerate}
\renewcommand{\theenumi}{\arabic{enumi}}
\renewcommand{\labelenumi}{($\dagger$\theenumi)}
\item
(Proposition 7.1) 
$\pi_{\lambda,\dagger}{\cong}\pi_{\lambda^*,\dagger}{\cong}
\pi_{\lambda',\dagger}{\cong}\pi_{\lambda'^*,\dagger}$. 
\item
(Proposition 7.3) 
\begin{enumerate}
\rm\item\it
If $\lambda^*=\lambda$, $\lambda'{\;\neq\;}\lambda$, $\lambda'^*{\;\neq\;}\lambda$, 
then $\pi_{\lambda,\sharp,\dagger}^+{\cong}\pi_{\lambda',\sharp,\dagger}^+$ and 
$\pi_{\lambda,\sharp,\dagger}^-{\cong}\pi_{\lambda',\sharp,\dagger}^-$. 
They are irreducible and mutually inequivalent. 
\rm\item\it
If $\lambda^*{\;\neq\;}\lambda$, $\lambda'=\lambda$, $\lambda'^*{\;\neq\;}\lambda$, 
then $\pi_{\lambda,\flat,\dagger}^+{\cong}\pi_{\lambda',\flat,\dagger}^+$ and 
$\pi_{\lambda,\flat,\dagger}^-{\cong}\pi_{\lambda',\flat,\dagger}^-$. 
They are irreducible and mutually inequivalent. 
\rm\item\it
If $\lambda^*{\;\neq\;}\lambda$, $\lambda'{\;\neq\;}\lambda$, $\lambda'^*=\lambda$, 
then $\pi_{\lambda,\natural,\dagger}^+{\cong}\pi_{\lambda',\natural,\dagger}^+$ and 
$\pi_{\lambda,\natural,\dagger}^-{\cong}\pi_{\lambda',\natural,\dagger}^-$. 
They are irreducible and mutually inequivalent. 
\end{enumerate}
\item
(Corollary 7.7) 
If $\lambda^*=\lambda'=\lambda'^*=\lambda$, then 
$\pi_{\lambda,\dagger}$ decomposes into four inequivalent 
subrepresentations $\pi_{\lambda,\dagger}^{++}$, $\pi_{\lambda,\dagger}^{+-}$, 
$\pi_{\lambda,\dagger}^{-+}$ and $\pi_{\lambda,\dagger}^{--}$ of the same degree. 
\item
(Theorem 7.9) 
Let $\{\lambda_i,\lambda_i^*,\lambda_i',\lambda_i'^*,
\mu_j,\mu_j{'},\nu_k,\nu_k^*,
\kappa_l,\kappa_l',
\iota_m\}_{i,j,k,l,m}$ 
be the set of $2$-tuples of Young diagrams 
such that 
\begin{equation*}
\begin{split}
&\lambda_i{'}{\neq}\lambda_i,\lambda_i^*{\neq}\lambda_i,\lambda_i'^*{\neq}\lambda_i,{\quad}
\mu_j^*=\mu_j,\mu_j{'}{\neq}\mu_j,\mu_j'^*{\neq}\mu_j,\\
&\nu_k^*{\neq}\nu_k,\nu_k{'}=\nu_k,\nu_k'^*{\neq}\nu_k,{\quad}
\kappa_l^*{\neq}\kappa_l,\kappa_l{'}{\neq}\kappa_l,\kappa_l'^*=\kappa_l,{\quad}
\iota_m{'}=\iota_m^*=\iota_m'^*=\iota_m.
\end{split}
\end{equation*}
Then 
\[
\{\pi_{\lambda_i,\dagger},\pi_{\mu_j,\sharp,\dagger}^{\pm},
\pi_{\nu_k,\flat,\dagger}^{\pm},
\pi_{\kappa_l,\natural,\dagger}^{\pm},
\pi_{\iota_m,\dagger}^{++},\pi_{\iota_m,\dagger}^{+-},\,
\pi_{\iota_m,\dagger}^{-+},\,\pi_{\iota_m,\dagger}^{--}\}_{i,j,k,l,m}
\]
is a basic set of irreducible representations of the intersection subalgebra. 
\end{enumerate}
\section{Preliminaries on group graded algebras}
In this section, we review some properties about group graded algebras. 
Statements listed below and its proof can be found in Curtis-Reiner \cite{C-R1}. 
For a detailed account of group graded rings, we refer to \cite{N-O}. 
Here we give the definition of the $S$-graded Clifford system in the algebra 
over the commutative ring for the finite group $S$ and give properties of 
$S$-graded Clifford systems. 
The definition of the $S$-graded $R$-algebra given in \cite{C-R1} is 
somewhat unpopular, so we adopt the term {\lq\lq}$S$-crossed product{\rq\rq} instead of 
{\lq\lq}$S$-graded $R$-algebra{\rq\rq} in accordance with \cite{N-O}. 
\begin{Def}
Let $S$ be a finite group, $R$ a commutative ring, and $A$ an $R$-algebra, finitely 
generated over $R$ as module. A family of $R$-submodules $\{A_s\}_{s{\in}S}$ of $A$, 
indexed by the elements of $S$, is called an \it $S$-graded Clifford system \rm in $A$, if 
the following conditions are satisfied.
\renewcommand{\theenumi}{\arabic{enumi}}
\renewcommand{\labelenumi}{(C\theenumi)}
\begin{enumerate}
\item\quad
$A_sA_t=A_{st}$
\item\quad
For each $s{\in}S$, there exists a unit $a_s{\in}A$ such that 
$A_s=a_sA_1=A_1a_s$
\item\quad
$A=\bigoplus_{s{\in}S}A_s$. 
\item\quad
$1{\in}A_1$
\end{enumerate}
\renewcommand{\theenumi}{\arabic{enumi}}
\renewcommand{\labelenumi}{(\theenumi)}
\end{Def}
An $R$-algebra satisfying these conditions is called an \it $S$-crossed product\rm. 
We notice that the $S$-graded Clifford system is a generalization of the 
twisted group algebra $R[S]_\alpha$ with factor set $\alpha$, and trivial 
$G$-action on $R$.
The next five results are basic facts about the representation theory of the 
$S$-graded Clifford system. We will use these results to derive the branching 
rule later. The proofs of these results are shown in Curtis-Reiner \cite{C-R1} Chapter 1, 
Section 11C (p.267-279).
\begin{Prop}[\cite{C-R1} (11.13)]
Let $A$ be an $R$-algebra with an $S$-graded Clifford system, and let $L$ be 
a left $A_1$-module, and $M$ a left $A$-module. Then there are two isomorphisms 
of $R$-modules 
\[ (1)\quad \operatorname{Hom}_{A_1}(L,M_{A_1})\cong \operatorname{Hom}_{A}(L^{A},M) \]
and
\[ (2)\quad \operatorname{Hom}_{A_1}(M_{A_1},L)\cong \operatorname{Hom}_{A}(M,L^{A}), \]
where $M_{A_1}$ is the $A_1$-module obtained by restriction of scalars from 
$A$ to $A_1$, and $L^A$ is the induced module from $A_1$ of $A$ defined as follows.
\[
L^A=A\otimes_{A_1}L
\]
\end{Prop}
A left $A_1$-module $L$ is called \it stable relative to $A$ \rm if $L$ is 
$A_1$-isomorphic to all $a_s\otimes L$ with $s\in S$. 
\begin{Prop}[\cite{C-R1} (11.14)]
Let $A$ have an $S$-graded Clifford system, let $L$ be a left $A_1$-module, and 
let $E$ denote the endomorphism algebra $\operatorname{End}_A(L^A)$, viewed as a 
ring of right operators on $L^A$. For each $s\in S$, let 
\[
E_s=\{f\in E \mid (1\otimes L)f\subseteq a_s\otimes L\}
\]
Then
\begin{enumerate}
\rm\item\it
For all $s,t\in S$, we have 
\[ A_s(a_t\otimes L)=a_{st}\otimes L, \quad (a_s\otimes L)E_s\subseteq (a_{st}\otimes L) \]
\[ E_sE_t\subseteq E_{st}, \quad 1\in E_1, \quad E=\oplus_{s\in S}E_s \]
\rm\item\it
Each element $\varphi\in \operatorname{Hom}_{A_1}(1\otimes L,a_s\otimes L)$ extends to a 
unique element $\Hat{\varphi}\in E_s$, given by 
$(a\otimes l)\Hat{\varphi}=a(1\otimes l)\varphi$ 
for $l\in L,a\in A$. The map $\varphi \rightarrow \Hat{\varphi}$ defines an isomorphism 
of $R$-modules
\[ \operatorname{Hom}_{A_1}(1\otimes L,a_s\otimes L)\cong E_s \quad s\in S \]
and this is an isomorphism of $R$-algebras when $a_s=1$.
\rm\item\it
If $L$ is stable and $R$ is noetherian, then $E$ has an $S$-graded Clifford system, 
with units $\Hat{e}_s\in E_s$ defined by (2) from $A_1$-isomorphisms 
$e_s:1\otimes L \rightarrow a_s\otimes L$, for all $s\in S$.
\end{enumerate}
\end{Prop}
\begin{Prop}[\cite{C-R1} (11.15)]
Let $A$ have an $S$-graded Clifford system over a field $K$, and let $L$ be a 
simple $A_1$-module. Then $E_1$ is a finite dimensional division algebra over $K$, 
and 
\[ T=\{t\in S \mid 1\otimes L \cong a_t\otimes L \quad \text{as $A_1$-modules}\} \]
is a subgroup of $S$. Further the $K$-algebra 
\[ E_T=\sum_{t\in T}E_t \]
is a $T$-crossed product.
\end{Prop}
The subgroup $T$ defined in Proposition 2.4 is called the \it inertial group \rm 
of $L$. 
\begin{Prop}[\cite{C-R1} (11.16)]
Let $A$ be an $S$-crossed product over a field $K$, $M$ a simple $A$-module, 
and $L$ a simple submodule of $M_{A_1}$. Then the following statements hold.
\begin{enumerate}
\rm\item\it
$M_{A_1}$ is a semisimple $A_1$-module, whose simple summands are isomorphic to 
the simple $A_1$-modules $\{a_s\otimes L \mid s \in S\}$.
\rm\item\it
Let $T=\{t\in S \mid L\cong a_t\otimes L \quad \text{as $A_1$-module}\}$, and let 
$B=\sum_{t\in T}A_t$. Then the homogeneous component $N$ of $M_{A_1}$ containing $L$ 
is a $B$-module, and $M$ is $A$-isomorphic to the induced module $A\otimes_B L$.
\rm\item\it
$N$ is a simple submodule of $L^B$, and $L$ is stable relative to $B$.
\end{enumerate}
\end{Prop}
\begin{Theo}[\cite{C-R1} (11.17)]
Let $A$ be an $S$-crossed product over a field $K$, and let $L$ be a simple $A_1$-module 
which is stable relative to $A$. 
Let $E=\operatorname{End}_A(L^A)$ as ring of right operators of $L^A$. 
\begin{enumerate}
\rm\item\it
There is an isomorphism $\theta$ from the lattice of left ideals $I$ in $E$ onto the 
lattice of $A$-submodules $U$ of $L^A$, given by
\[ U=L^A\cdot I,\quad I=\{\gamma \in E \mid L^A\gamma\subseteq U\} \]
\rm\item\it
The $A$-module $L^A\cdot I$ corresponding to $I$ is $A$-isomorphic to $L\otimes_{E_1}I$, 
where the action of $A$ on this tensor product is given by 
\[ a(l\otimes \gamma)=(al)\Hat{e}_s^{-1}\otimes \Hat{e}_s\gamma 
\quad a\in A_s, l\in L, \gamma\in I \]
with $E_1$ and the units $\{\Hat{e}_s\}$ in $E_s$ defined as in Proposition 2.3. 
Furthermore,
\[ \dim_K(L\otimes_{E_1}I)=(\dim_{E_1}L)(\dim_KI) \]
\rm\item\it
The lattice isomorphism $\theta$ is functorial, in the sense that \\
$E$-homomorphisms 
$f : I{\longrightarrow}I'$ between left ideals of $E$, correspond bijectively to 
$A$-homomorphisms 
$1{\otimes}f : L{\otimes}_{E_1}I{\longrightarrow}L{\otimes}_{E_1}I'$. 
\end{enumerate}
\end{Theo}
If $S=\mathbb{Z}_2=\{\Bar{0},\,\Bar{1}\}$ and $K$ is an algebraic closed field, then 
we obtain the following theorem which we will use later frequently. 
\begin{Theo}
Let $A$ be a $\mathbb{Z}_2$-crossed product over an algebraic closed field $K$, 
and $M$ a simple left $A$-module. Then one of the following holds : 
\begin{enumerate}
\rm\item\it
$M_{A_{\Bar{0}}}$ is a simple left $A_{\Bar{0}}$-module. The induced module 
$(M_{A_{\Bar{0}}})^A$ is a direct sum of two mutually non-isomorphic simple left 
$A$-modules $M^1$ and $M^2$. 
One of those is $A$-isomorphic to $M$ and their restrictions $M^1_{A_{\bar{0}}}$ and 
$M^2_{A_{\bar{0}}}$ are $A_1$-isomorphic to $M_{A_{\Bar{0}}}$. 
\rm\item\it
$M_{A_{\Bar{0}}}=L_1{\oplus}L_2$ where $L_1,L_2$ are mutually non-isomorphic 
simple left $A_{\Bar{0}}$-modules which have same dimension. 
The induced modules $L_1^A$ and $L_2^A$ are 
$A$-isomorphic to $M$. Both $L_1$ and $L_2$ never appear as irreducible constituents of 
$M'_{A_0}$, where $M'$ is any simple left $A$-module which is non-isomorphic to $M$. 
\end{enumerate}
Furthermore, all $M_{A_{\bar{0}}}$ in the case (1) and all $L_1$ and $L_2$ 
in the case (2) are mutually non-isomorphic except for the pairs 
$(M^1_{A_{\bar{0}}},M^2_{A_{\bar{0}}})$ where $M$ is in the case (1). 
\end{Theo}
\begin{proof}
Let $L$ be a simple submodule of $M_{A_{\Bar{0}}}$ and $T$ the inertial group of $L$. 
Since $T$ is a subgroup of $\mathbb{Z}_2$, $T$ must be $\{0\}$ or $\mathbb{Z}_2$. \par
If $T=\{0\}$, then $a_{\Bar{1}}{\otimes}L{\ncong}1{\otimes}L$ as 
$A_{\bar{0}}$-modules. 
From Proposition 2.5 (1) we have 
\[
M_{A_{\Bar{0}}}{\cong}m_{\Bar{0}}(1{\otimes}L)
{\oplus}m_{\Bar{1}}(a_{\Bar{1}}{\otimes}L)
\]
as $A_{\Bar{0}}$-modules where $m_{\Bar{0}}$ and 
$m_{\Bar{1}}$ are multiplicities of $1{\otimes}L$ and $a_{\Bar{1}}{\otimes}L$ respectively. 
We also have
\[
M_{A_{\Bar{0}}}=a_{\Bar{1}}M_{A_{\Bar{0}}}{\cong}m_{\Bar{0}}(a_{\Bar{1}}{\otimes}L)
{\oplus}m_{\Bar{1}}(1{\otimes}L)
\]
Thus $m_{\Bar{0}}=m_{\Bar{1}}$, so we may write 
$M_{A_{\Bar{0}}}{\cong}m\big{\{}(1{\otimes}L)
{\oplus}(a_{\Bar{1}}{\otimes}L)\big{\}}$ for some positive integer $m$. 
Hence $\dim_{K}M=2m{\dim_{K}}L$. 
On the other hand, $\dim_{K}M{\leq}2{\dim_{K}}L$ since $M$ is an irreducible constituent 
of $L^A=(1{\otimes}L){\oplus}(a_{\Bar{1}}{\otimes}L)$. Therefore $m=1$ and 
$L^A=M$ holds. 
Let $M'$ be any simple left $A$-module which is non-isomorphic to $M$. 
If $M'_{A_{\bar{0}}}$ contains $L_1$ 
or $L_2$ as an irreducible constituent, say $L_1$ is contained, 
then $(L_1)^A$ is $A$-isomorphic to both $M$ 
and $M'$, but this is contradiction. 
This corresponds to the case (2). \par
If $T=\mathbb{Z}_2$, then $a_{\Bar{1}}{\otimes}L{\cong}1{\otimes}L$ as 
$A_{\bar{0}}$-modules and $M_{A_{\Bar{0}}}{\cong}mL$ for some positive integer $m$. 
$E=\operatorname{End}_{A}L^A$ has a $\mathbb{Z}_2$-graded Clifford system by Proposition 2.3 (3). 
Hence we may write 
\[
E=E_{\bar{0}}{\otimes}E_{\bar{1}}=E_{\bar{0}}{\otimes}E_{\bar{0}}\hat{e}_{\bar{1}}
\]
Moreover, since $K$ is algebraic closed we have $E_{\bar{0}}{\cong}K$. 
Multiplying $\hat{e}_{\bar{1}}$ by a suitable scalar, we may assume $\hat{e}_{\bar{1}}^2=1$. 
Thus $E$ is isomorphic to the quotient algebra $K[X]/(X^2-1)$ of the polynomial ring 
by the correspondence $\hat{e}_{\bar{1}}{\mapsto}X$. 
Let $I_1,I_2$ be ideals of $E$ which are defined by 
\[
I_1={\langle}1+\hat{e}_{\bar{1}}{\rangle},{\quad}
I_2={\langle}1-\hat{e}_{\bar{1}}{\rangle}.
\]
Then we readily see that $\dim_{K}I_1=\dim_{K}I_2=1$ and $E=I_1{\oplus}I_2$. 
Since $I_1$ and $I_2$ are mutually non-isomorphic simple $E$-modules, 
$L^A{\cdot}I_1$ and $L^A{\cdot}I_2$ are mutually non-isomorphic simple $A$-modules 
by Theorem 2.6 (1),(3). 
From Theorem 2.6 (2), we have 
\begin{equation*}
\begin{split}
\dim_{K}L{\cdot}I_1&=\dim_{K}L{\times}\dim_{K}I_1=\dim_{K}L,\\
\dim_{K}L{\cdot}I_2&=\dim_{K}L{\times}\dim_{K}I_2=\dim_{K}L.
\end{split}
\end{equation*}
Since $\dim_{K}L^A=2\dim_{K}L$, $L^A$ decomposes into two simple submodules as follows. 
\[
L^A=L{\cdot}I_1{\oplus}L{\cdot}I_2
\]
$M$ is an irreducible constituent of $L^A$, hence $M{\cong}L{\cdot}I_1$ or 
$M{\cong}L{\cdot}I_2$. 
In each case, the restriction $M_{A_{\bar{0}}}$ is isomorphic to $L$ by Theorem 2.6 (2). 
This corresponds to the case (1). \par
For the last statement, we have already shown that 
$M^1_{A_{\bar{0}}}{\cong}M^2_{A_{\bar{0}}}$ and $L_1{\ncong}L_2$. 
Other pairs of simple $A_1$-modules are mutually non-isomorphic 
since their induced modules are mutually non-isomorphic. 
Thus we have completed the proof. 
\end{proof}
\section{Hecke algebras of type $B_n$ and $D_n$}
Let $(W,S=\{s_1,\hdots,s_n\})$ be a Coxeter system of rank $n$. 
Let $R$ be a commutative domain with $1$, and let $q_i(i=1,\hdots,n)$ be any invertible 
elements of $R$ such that $q_i=q_j$ if $s_i$ is conjugate to $s_j$ in $W$.
The Iwahori-Hecke algebra $\mathscr{H}_{R}(W,S)$ is an $R$-algebra generated 
by $\{T_{s_i}|s_i{\in}S\}$ with the relations: 
\renewcommand{\theenumi}{\arabic{enumi}}
\renewcommand{\labelenumi}{(H\theenumi)}
\begin{enumerate}
\item
$T_{s_i}^2 = (q_i-1)T_{s_i}+q_i$ \qquad if $i=1,2,\hdots,n$,
\item
$(T_{s_i}T_{s_j})^{k_{ij}}=(T_{s_j}T_{s_i})^{k_{ij}}$ \qquad if $m_{ij}=2k_{ij}$,
\item
$(T_{s_i}T_{s_j})^{k_{ij}}T_{s_i}=(T_{s_j}T_{s_i})^{k_{ij}}T_{s_j}$ \qquad if $m_{ij}=2k_{ij}+1$,
\end{enumerate}
\renewcommand{\theenumi}{\arabic{enumi}}
\renewcommand{\labelenumi}{(\theenumi)}
where $m_{ij}$ is the order of $s_is_j$ in $W$. 
We define $T_w=T_{s_{i_1}}T_{s_{i_2}}{\cdots}T_{s_{i_k}}$ where 
$w=s_{i_1}s_{i_2}{\cdots}s_{i_k}$ is 
a reduced expression of $w$. It is known that $\{T_w|w{\in}W\}$ form a basis of $\mathscr{H}_{R}(W,S)$ as free $R$-modules. 
The relations (H1)--(H3) is equivalent to the following two relations: 
\renewcommand{\theenumi}{\arabic{enumi}}
\renewcommand{\labelenumi}{(h\theenumi)}
\begin{enumerate}
\item
$T_{s_i}T_w = T_{s_iw}$ \qquad if $l(w)<l(s_iw)$,
\item
$T_{s_i}T_w = (q_i-1)T_w + q_iT_{s_iw}$ \qquad if $l(w)>l(s_iw)$,
\end{enumerate}
or equivalently, 
\renewcommand{\theenumi}{\arabic{enumi}}
\renewcommand{\labelenumi}{(h'\theenumi)}
\begin{enumerate}
\item
$T_{w}T_{s_i} = T_{ws_i}$ {\qquad} if $l(w)<l(ws_i)$,
\item
$T_wT_{s_i} = (q_i-1)T_w + q_iT_{ws_i}$ {\qquad} if $l(w)>l(ws_i)$,
\end{enumerate}
\renewcommand{\theenumi}{\arabic{enumi}}
\renewcommand{\labelenumi}{(\theenumi)}
where $l(w)$ means the length of $w$. \par
We set $u=q_1$, $q=q_2=q_3={\cdots}q_n$. 
The Hecke algebra $\mathscr{H}_{R,B_n}(u,q)$ of type $B_n$ $(n{\geq}2)$ is the algebra 
over $R$ defined by the generators $a_1,a_2,\hdots,a_n$ and the following defining relations.
\renewcommand{\theenumi}{\arabic{enumi}}
\renewcommand{\labelenumi}{(B\theenumi)}
\begin{enumerate}
\item
$a_1^2 = (u-1)a_1+u$
\item
$a_i^2 = (q-1)a_i+q \qquad$ if $i=2,3,\hdots,n$
\item
$a_1a_2a_1a_2 = a_2a_1a_2a_1$
\item
$a_ia_{i+1}a_i = a_{i+1}a_ia_{i+1} \qquad$ if $i=2,3,\hdots,n-1$
\item
$a_ia_j = a_ja_i \qquad$ if $|i-j|>1$
\end{enumerate}
\renewcommand{\theenumi}{\arabic{enumi}}
\renewcommand{\labelenumi}{(\theenumi)}
Let $u=1$ and consider $\mathscr{H}_{R,B_n}(1,q)$. Then (B1) is reduced to (B1') $a_1^2=1$.
One can readily check that the element $\Bar{a}_1=a_1a_2a_1$ satisfies $\Bar{a}_1^2=(q-1)\Bar{a}_1+q$. 
For $n{\geq}2$, we define $\Bar{\mathscr{H}}_{R,B_n}(1,q)$ to be the subalgebra of 
$\mathscr{H}_{R,B_n}(1,q)$ generated by $\Bar{a}_1,\Bar{a}_2,{\hdots},\Bar{a}_n$ 
where $\Bar{a}_i=a_i$ for $i>1$. 
We also define $\Bar{\mathscr{H}}_{R,B_1}(1,q)$ to be the subalgebra of 
$\mathscr{H}_{R,B_1}(1,q)$ generated by the identity element. 
We readily see that $\Bar{\mathscr{H}}_{R,B_n}(1,q)$ satisfies the following relations. 
\renewcommand{\theenumi}{\arabic{enumi}}
\renewcommand{\labelenumi}{(D\theenumi)}
\begin{enumerate}
\item
$\Bar{a}_i^2 = (q-1)\Bar{a}_i+q \qquad$ if $i=1,2,\hdots,n$
\item
$\Bar{a}_1\Bar{a}_i=\Bar{a}_i\Bar{a}_1 \qquad$ if $i{\neq}3$
\item
$\Bar{a}_1\Bar{a}_3\Bar{a}_1 = \Bar{a}_3\Bar{a}_1\Bar{a}_3$
\item
$\Bar{a}_i\Bar{a}_{i+1}\Bar{a}_i = \Bar{a}_{i+1}\Bar{a}_i\Bar{a}_{i+1} \qquad$ 
if $i=2,3,\hdots,n-1$
\item
$\Bar{a}_i\Bar{a}_j = \Bar{a}_j\Bar{a}_i \qquad$ if $|i-j|>1$
\end{enumerate}
\renewcommand{\theenumi}{\arabic{enumi}}
\renewcommand{\labelenumi}{(\theenumi)}
For $n{\geq}4$, it is known that (D1)--(D5) are defining relations of the Hecke algebra $\mathscr{H}_{R,D_n}(q)$ of type $D_n$, 
so we may identify $\mathscr{H}_{R,D_n}(q)$ with $\Bar{\mathscr{H}}_{R,B_n}(1,q)$. \par
Let $u$ and $q$ be indeterminates and 
\[ R_0=\mathbb{Z}[u^{\pm{1}},q^{\pm{1}}],{\quad}R_1=\mathbb{Z}[q^{\pm{1}}]. \]
It is known that $\mathscr{H}_{R_0,B_n}(u,q)$ is a free $R_0$-module of rank $2^nn!$ 
\cite{Ariki-Koike,Hoefsmit} and that $\Bar{\mathscr{H}}_{R_1,B_n}(1,q)$ is a free $R_1$-module of 
rank $2^{n-1}n!$. Accordingly, $\mathscr{H}_{R_1,B_n}(1,q)$ is also a free $R_1$-module of rank $2^nn!$. 
\section{Fixed subalgebras of $\mathscr{H}_{R_3,B_n}(1,q)$ by involutions}
It is obvious that there is an algebra automorphism $\sharp$ of $\mathscr{H}_{R_1,B_n}(1,q)$ 
of order $2$ which is defined by $a_1^{\sharp}=-a_1$ and $a_i^{\sharp}=a_i$ ($i{\geq}2$). 
We define $\mathscr{H}_{R_1,B_n}(1,q)^{\sharp}$ to be the subalgebra consisting of fixed elements 
of $\sharp$ as follows. 
\[
\mathscr{H}_{R_1,B_n}(1,q)^{\sharp}=\{h{\in}\mathscr{H}_{R_1,B_n}(1,q)\,|\,h^{\sharp}=h\}.
\]
The relations (B1') and (B2)--(B5) imply that every monomial with occurrences of even 
(resp. odd) numbers of $a_1$ also has occurrences of even (resp. odd) numbers of $a_1$ 
in each term of any other expression of it. 
Therefore $\mathscr{H}_{R_1,B_n}(1,q)^{\sharp}$ is generated by monomials with occurrences 
of even numbers of $a_1$. 
Since the Weyl group $W_{D_n}$ of type $D_n$ is a normal subgroup of the Weyl group $W_{B_n}$ 
of type $B_n$ of index $2$, and consists of all elements which can be written as 
products of $s_1,s_2,{\hdots},s_n$ with even numbers of factors of $s_1$. 
On the other hand, $W_{D_n}$ is generated by $s_1s_2s_1,s_2,s_3,{\hdots},s_n$ 
(see for example \cite{G-P},1.4.8). 
Thus we have the following direct sum of $R_1$-modules: 
\[
\mathscr{H}_{R_1,B_n}(1,q)^{\sharp}=
\bigoplus_{w{\in}W_{D_n}}R_1T_{w}=\Bar{\mathscr{H}}_{R_1,B_n}(1,q).
\]
The argument above is valid even if $n=2,3$. 
$\bigoplus_{w{\in}W_{D_n}}RT_{w}$ is closed under the product by (B1') and (B2)--(B5), 
Hence $\Bar{\mathscr{H}}_{R_1,B_n}(1,q)=\mathscr{H}_{R_1,B_n}(1,q)^{\sharp}$ holds.
Let $R_2,R_3$ be the rings defined by 
\begin{equation*}
\begin{split}
R_2&=\mathbb{Z}[u^{\pm{1}},q^{\pm{1}},(u+1)^{-1},(q+1)^{-1},\dfrac{1}{2}]\\
R_3&=\mathbb{Z}[q^{\pm{1}},(q+1)^{-1},\dfrac{1}{2}].\\
\end{split}
\end{equation*}
It is known that there exists an algebra automorphism $\natural$ of 
$\mathscr{H}_{R_0,B_n}(u,q)$ of order $2$ which is defined by 
$a_1^{\natural}=(u-1)-a_1$ and $a_i^{\natural}=(q-1)-a_i$ ($i{\geq}2$). 
$\natural$ is called Goldman involution. 
We consider the fixed subalgebra by $\natural$ over $R_2$. 
We define the elements $b_i$ ($i=1,2,\hdots,n$) of 
$\mathscr{H}_{R_2,B_n}(u,q)$ by 
\[
b_1=\dfrac{a_1-a_1^{\natural}}{u+1},{\qquad}
b_i=\dfrac{a_i-a_i^{\natural}}{q+1}{\qquad}\text{if $i=2,3,\hdots,n$}
\]
Then the following holds.
\begin{Prop}[\cite{Mitsuhashi2}, Proposition 3.2]
$b_1,b_2,{\hdots}b_n$ generate 
$\mathscr{H}_{R_2,B_n}(u,q)$ and constitute with the relations
\renewcommand{\theenumi}{\arabic{enumi}}
\renewcommand{\labelenumi}{(B'\theenumi)}
\begin{enumerate}
\rm \item \it
$b_i^2=1{\qquad}$ if $i=1,2,\hdots,n$
\rm \item \it
$b_1b_2b_1b_2=b_2b_1b_2b_1-2\dfrac{(u-1)(q-1)}{(u+1)(q+1)}(b_1b_2-b_2b_1)$
\rm \item \it
$b_ib_{i+1}b_i=b_{i+1}b_ib_{i+1}-\Big{(}\dfrac{q-1}{q+1}\Big{)}^2(b_i-b_{i+1})\qquad$ 
if $i=2,3,\hdots,n$
\rm \item \it
$b_ib_j=b_jb_i\qquad$ if $|i-j|>1$
\end{enumerate}
\renewcommand{\theenumi}{\arabic{enumi}}
\renewcommand{\labelenumi}{(\theenumi)}
a presentation of $\mathscr{H}_{R_2,B_n}(u,q)$.
\end{Prop}
Consider the following sets of monomials.
\begin{equation*}
\begin{split}
\mathcal{S}_1 &= \{1,b_1\} \\
\mathcal{S}_2 &= \{1,b_2,b_2b_1,b_2b_1b_2\} \\
\vdots \\
\mathcal{S}_i &= \{1,b_i,b_ib_{i-1},\hdots ,b_ib_{i-1}\cdots b_2b_1,b_ib_{i-1}\cdots b_2b_1b_2,
\hdots,\\ 
&{\hspace{2cm}}b_ib_{i-1}\cdots b_2b_1b_2\cdots b_{i-1}b_i\} \\
\vdots \\
\mathcal{S}_n &= \{1,b_n,b_nb_{n-1},\hdots ,b_nb_{n-1}\cdots b_2b_1,
b_nb_{n-1}\cdots b_2b_1b_2,
\hdots,\\
&{\hspace{2cm}}b_nb_{n-1}\cdots b_2b_1b_2\cdots b_{n-1}b_n\}
\end{split}
\end{equation*}
We shall say that $M=U_1U_2\cdots U_{n}$ is a
\it monomial in $b_i$-normal form \rm in $\mathscr{H}_{R_2,B_n}(u,q)$, 
if $U_i \in \mathcal{S}_i$ for $i=1,2,\hdots ,n$. 
Then we have  
\begin{Prop}[\cite{Mitsuhashi2},Proposition 3.3]
\[ \mathscr{H}_{R_2,B_n}(u,q) = \bigoplus_{U_i{\in}\mathcal{S}_i}R_2U_1U_2{\cdots}U_n \]
\end{Prop}
Applying the case $u=1$ to the above proposition, we also obtain 
\[
\mathscr{H}_{R_3,B_n}(1,q) = \bigoplus_{U_i{\in}\mathcal{S}_i}R_3U_1U_2{\cdots}U_n.
\]
One can determine the fixed subalgebra 
by $\natural$ in $\mathscr{H}_{R_2,B_n}(u,q)$. 
Let $\mathscr{E}_n^{\natural}$ be the set of all monomials in $b_i$-normal form in 
$\mathscr{H}_{R_2,B_n}(u,q)$ which are products of even numbers of $b_i$'s. 
We define $\mathscr{H}_{R_2,B_n}(u,q)^{\natural}$ to be the 
$R_2$-submodule of $\mathscr{H}_{R_2,B_n}(u,q)$ defined as follows.
\[ 
\mathscr{H}_{R_2,B_n}(u,q)^{\natural}=\bigoplus_{M\in\mathscr{E}_n^{\natural}}R_2M 
\]
\begin{Prop}[\cite{Mitsuhashi2},Proposition 3.6]
$\mathscr{H}_{R_2,B_n}(u,q)^{\natural}$ is the subalgebra of $\mathscr{H}_{R_2,B_n}(u,q)$ 
which consists of all the elements fixed under $\natural$. 
Furthermore $\operatorname{rank}_{R_2}\mathscr{H}_{R_2,B_n}(u,q)^{\natural}=2^{n-1}n!$.
\end{Prop}
We mention that this proposition may be applied to the case $u=1$, and yields 
the same assertion for $\mathscr{H}_{R_3,B_n}(1,q)$ and 
$\mathscr{H}_{R_3,B_n}(1,q)^{\natural}$. In this case, 
$\natural$ turns out to be the automorphism determined by 
$a_1^{\natural}=-a_1$ and $a_i^{\natural}=(q-1)-a_i$ ($i{\geq}2$). 
Let $\mathscr{S}_n$ be the set of all monomials in $b_i$-normal form in 
$\mathscr{H}_{R_2,B_n}(u,q)$. 
We define $\mathscr{H}_{R_2,B_n}(u,q)^{-\natural}$ to be the $R_2$-submodule of 
$\mathscr{H}_{R_2,B_n}(u,q)$ as follows.
\[ \mathscr{H}_{R_2,B_n}(u,q)^{-\natural}
=\bigoplus_{M\in\mathscr{S}_n\backslash\mathscr{E}_n^{\natural}}R_2M \]
Then the following statement holds.
\begin{Prop}[\cite{Mitsuhashi2}, Proposition 3.7]
$\mathscr{H}_{R_2,B_n}(u,q)$ is a 
$\mathbb{Z}_{2}$-crossed product with $R_2$-submodules 
$A_{\Bar{0}}=\mathscr{H}_{R_2,B_n}(u,q)^{\natural}$, 
$A_{\Bar{1}}=\mathscr{H}_{R_2,B_n}(u,q)^{-\natural}$ and units $1{\in}A_{\Bar{0}}$, 
$b_1{\in}A_{\Bar{1}}$.
\end{Prop}
If we take $u=1$, we obtain a $\mathbb{Z}_{2}$-crossed product 
\[ 
\mathscr{H}_{R_3,B_n}(1,q)
=\mathscr{H}_{R_3,B_n}(1,q)^{\natural}{\oplus}\mathscr{H}_{R_3,B_n}(1,q)^{-\natural}
\]
Besides $\sharp$, $\natural$, 
we define one more algebra automorphism of $\mathscr{H}_{R_3,B_n}(1,q)$ of order $2$. 
Let $\flat$ be the map defined by $a_1^{\flat}=a_1$ and 
$a_i^{\flat}=(q-1)-a_i$ ($i{\geq}2$). Then $\flat$ can be extended to an 
algebra automorphism of $\mathscr{H}_{R_1,B_n}(1,q)$. 
One can readily see that $\sharp{\cdot}\flat=\flat{\cdot}\sharp=\natural$. 
We consider the fixed subalgebra by $\flat$ over $R_3$. 
We define $\mathscr{H}_{R_3,B_n}(1,q)^{\flat}$ to be the subalgebra consisting of fixed elements 
of $\flat$ as follows. 
\[
\mathscr{H}_{R_3,B_n}(1,q)^{\flat}=\{h{\in}\mathscr{H}_{R_3,B_n}(1,q)\,|\,h^{\flat}=h\}.
\]
Let $\mathscr{E}_n^{\flat}$ be the set of all monomials in $b_i$-normal form in 
$\mathscr{H}_{R_3,B_n}(1,q)^{\flat}$ which have occurrences 
of even numbers of $b_i$ ($i{\neq}1$). 
For $n{\geq}1$, we define $\Tilde{\mathscr{H}}_{R_3,B_n}(1,q)$ to be the 
$R_3$-submodule of $\mathscr{H}_{R_3,B_n}(1,q)$ defined as follows.
\[ 
\Tilde{\mathscr{H}}_{R_3,B_n}(1,q)=\bigoplus_{M\in\mathscr{E}_n^{\flat}}R_3M 
\]
We have $b_1^{\flat}=b_1$ and $b_i^{\flat}=-b_i$ ($i{\geq}2$) immediately. 
Hence the following holds. 
\begin{Prop}
$\mathscr{H}_{R_3,B_n}(1,q)^{\flat}=\Tilde{\mathscr{H}}_{R_3,B_n}(1,q)$. Furthermore, \\
$\operatorname{rank}_{R_3}\mathscr{H}_{R_3,B_n}(1,q)^{\flat}=2^{n-1}n!$.
\end{Prop}
\begin{proof}
Equality of two algebras is by Proposition 4.1. 
We notice that $|\mathcal{S}_i|=2i$. Let $\mathcal{S}^o_i$ be the subset of 
$\mathcal{S}_i$ consisting of the monomials with occurrences of odd numbers of $b_i$ 
($i{\neq}1$). 
$|\mathcal{S}^o_i|$ is just $i$ if $i>1$ and $|\mathcal{S}^o_1|=0$. 
A monomial $U_1U_2{\cdots}U_n$ in $b_i$-normal form has occurrences of even numbers 
of $b_i$ ($i{\neq}1$) if and only if the number of $U_i$ such that 
$U_i\in\mathcal{S}^o_i$ is even. 
$U_1$ never belongs to $\mathcal{S}^o_1$. Hence there are $2^{n-2}$ cases 
of being so. In each case, there are $2n!$ ways of taking elements from $\mathcal{S}_i$'s. 
Thus, we have that there exist $2^{n-1}n!$ monomials in $b_i$-normal form with occurrences 
of even numbers of $b_i$ ($i{\neq}1$). 
\end{proof}
We define $\mathscr{H}_{R_3,B_n}(1,q)^{-\flat}$ to be the $R_3$-submodule of 
$\mathscr{H}_{R_3,B_n}(1,q)$ as follows. 
\[
\mathscr{H}_{R_3,B_n}(1,q)^{\flat}=\{h{\in}\mathscr{H}_{R_3,B_n}(1,q)\,|\,h^{\flat}=-h\}.
\]
Then by Proposition 4.1, we have 
\[
\mathscr{H}_{R_3,B_n}(1,q)^{-\flat}=\bigoplus_{M\in\mathscr{S}_n\backslash\mathscr{E}_n^{\flat}}R_3M.
\]
\begin{Prop}
$\mathscr{H}_{R_3,B_n}(1,q)$ is a 
$\mathbb{Z}_{2}$-crossed product with $R_3$-submodules 
$A_{\Bar{0}}=\mathscr{H}_{R_3,B_n}(1,q)^{\flat}$, 
$A_{\Bar{1}}=\mathscr{H}_{R_3,B_n}(1,q)^{-\flat}$ and units $1{\in}A_{\Bar{0}}$, 
$b_2{\in}A_{\Bar{1}}$.
\end{Prop}
\begin{proof}
Clearly 
\[
A_{\Bar{0}}A_{\Bar{0}}=A_{\Bar{0}},{\quad}
A_{\Bar{1}}A_{\Bar{1}}=A_{\Bar{0}},{\quad}
A_{\Bar{0}}A_{\Bar{1}}=A_{\Bar{1}}A_{\Bar{0}}=A_{\Bar{1}}
\]
hold. Hence (C1) is satisfied. 
We may take units $1{\in}A_{\Bar{0}}$ and $b_2{\in}A_{\Bar{1}}$ so that 
$b_2A_{\Bar{0}}=A_{\Bar{0}}b_2=A_{\Bar{1}}$. Therefore (C2) and (C4) hold. 
Since 
$\mathscr{H}_{R_3,B_n}(1,q)^{\flat}=\bigoplus_{M\in\mathscr{E}_n^{\flat}}R_3M$ and 
$\mathscr{H}_{R_3,B_n}(1,q)^{-\flat}=\bigoplus_{M\in\mathscr{S}_n\backslash\mathscr{E}_n^{\flat}}R_3M$, 
(C3) holds. 
\end{proof}
Now we apply the generators $b_i$ ($i=1,{\hdots},n$) to $\mathscr{H}_{R_3,B_n}(1,q)^{\sharp}$. 
We immediately have $b_1=a_1$ and $b_1^{\sharp}=-b_1,\,b_i^{\sharp}=b_i$ ($i=2,{\hdots},n$). 
Let $\mathscr{E}_n^{\sharp}$ be the set of all monomials in $b_i$-normal form in 
$\mathscr{H}_{R_3,B_n}(1,q)$ which have 
occurrences of even numbers of $b_1$. Then by Proposition 4.1, we immediately have 
\begin{equation*}
\begin{split}
\mathscr{H}_{R_3,B_n}(1,q)^{\sharp}&=\bigoplus_{M\in\mathscr{E}_n^{\sharp}}R_3M \\
\mathscr{H}_{R_3,B_n}(1,q)^{-\sharp}
&=\bigoplus_{M\in\mathscr{S}_n\backslash\mathscr{E}_n^{\sharp}}R_3M 
\end{split}
\end{equation*}
for $n{\geq}1$. Moreover we obtain
\begin{Prop}
$\mathscr{H}_{R_3,B_n}(1,q)$ is a 
$\mathbb{Z}_{2}$-crossed product with $R_3$-submodules 
$A_{\Bar{0}}=\mathscr{H}_{R_3,B_n}(1,q)^{\sharp}$, 
$A_{\Bar{1}}=\mathscr{H}_{R_3,B_n}(1,q)^{-\sharp}$ and units $1{\in}A_{\Bar{0}}$, 
$b_1{\in}A_{\Bar{1}}$.
\end{Prop}
\begin{proof}
One can obtain the proof in the same manner as in Proposition 4.6. 
\end{proof}
\section{A fixed subalgebra of $\mathscr{H}_{R_3,B_n}(1,q)^{\sharp}$ by Goldman involution}
Let $\Bar{b}_1=b_1b_2b_1,\,\Bar{b}_i=b_i$ ($i>1$). 
By a direct computation, we get the following. 
\begin{Prop}
$\Bar{b}_1,\Bar{b}_2,\hdots,\Bar{b}_n$ generate 
$\mathscr{H}_{R_3,B_n}(1,q)^{\sharp}$ and constitute with the relations
\renewcommand{\theenumi}{\arabic{enumi}}
\renewcommand{\labelenumi}{(D'\theenumi)}
\begin{enumerate}
\rm \item \it
$\Bar{b}_i^2=1{\qquad}$ if $i=1,2,\hdots,n$
\rm \item \it
$\Bar{b}_1\Bar{b}_i=\Bar{b}_i\Bar{b}_1 \qquad$ if $i{\neq}3$
\rm \item \it
$\Bar{b}_1\Bar{b}_3\Bar{b}_1=\Bar{b}_3\Bar{b}_1\Bar{b}_3
-\Big{(}\dfrac{q-1}{q+1}\Big{)}^2(\Bar{b}_1-\Bar{b}_3)$
\rm \item \it
$\Bar{b}_i\Bar{b}_{i+1}\Bar{b}_i=\Bar{b}_{i+1}\Bar{b}_i\Bar{b}_{i+1}
-\Big{(}\dfrac{q-1}{q+1}\Big{)}^2(\Bar{b}_i-\Bar{b}_{i+1})\qquad$ 
if $i=2,3,\hdots,{n-1}$
\rm \item \it
$\Bar{b}_i\Bar{b}_j=\Bar{b}_j\Bar{b}_i\qquad$ if $|i-j|>1$
\end{enumerate}
\renewcommand{\theenumi}{\arabic{enumi}}
\renewcommand{\labelenumi}{(\theenumi)}
a presentation of $\mathscr{H}_{R_3,B_n}(1,q)^{\sharp}$.
\end{Prop}
The relations (D'1)--(D'5) imply that every product of an even 
(resp. odd) number of generators $\Bar{b}_1,\Bar{b}_2,{\hdots},\Bar{b}_n$ 
must be written as a linear combination of products of even (resp. odd) numbers of the 
generators in any other expression of it. 
Therefore we define $\mathscr{H}_{R_3,B_n}(1,q)^{\dagger}$ 
to be the $R_3$-subalgebra of $\mathscr{H}_{R_3,B_n}(1,q)^{\sharp}$ 
generated by all the monomials of even numbers of factors of 
$\Bar{b}_1,\Bar{b}_2,{\hdots},\Bar{b}_n$. 
Goldman involution of $\mathscr{H}_{R_3,B_n}(1,q)^{\sharp}$ is given by 
$\Bar{a}_i^{\natural}=(q-1)-\Bar{a}_i$ ($i=1,2,{\hdots},n$). 
We readily see that $\Bar{b}_i^{\natural}=-\Bar{b}_i$. 
Thus the fixed subalgebra $(\mathscr{H}_{R_3,B_n}(1,q)^{\sharp})^{\natural}$ of 
$\natural$ coincides with $\mathscr{H}_{R_3,B_n}(1,q)^{\dagger}$. \par
We consider intersections of these subalgebras 
$\mathscr{H}_{R_3,B_n}(1,q)^{\natural},\mathscr{H}_{R_3,B_n}(1,q)^{\flat}$ and 
$\mathscr{H}_{R_3,B_n}(1,q)^{\sharp}$. 
We immediately have 
\[
\mathscr{H}_{R_3,B_n}(1,q)^{\dagger}=(\mathscr{H}_{R_3,B_n}(1,q)^{\sharp})^{\natural}
=\mathscr{H}_{R_3,B_n}(1,q)^{\sharp}{\cap}\mathscr{H}_{R_3,B_n}(1,q)^{\natural}.
\]
\begin{Prop}
\begin{equation*}
\begin{split}
\mathscr{H}_{R_3,B_n}(1,q)^{\natural}{\cap}\mathscr{H}_{R_3,B_n}(1,q)^{\flat}
&=\mathscr{H}_{R_3,B_n}(1,q)^{\flat}{\cap}\mathscr{H}_{R_3,B_n}(1,q)^{\sharp}\\
&=\mathscr{H}_{R_3,B_n}(1,q)^{\sharp}{\cap}\mathscr{H}_{R_3,B_n}(1,q)^{\natural}
\end{split}
\end{equation*}
\end{Prop}
\begin{proof}
Since $\flat=\natural{\cdot}\sharp$, we have
\begin{equation*}
\begin{split}
\mathscr{H}_{R_3,B_n}(1,q)^{\natural}{\cap}\mathscr{H}_{R_3,B_n}(1,q)^{\flat}
&=\mathscr{H}_{R_3,B_n}(1,q)^{\natural}{\cap}\mathscr{H}_{R_3,B_n}(1,q)^{\natural{\cdot}\sharp}\\
&{\subseteq}\mathscr{H}_{R_3,B_n}(1,q)^{\natural}{\cap}\mathscr{H}_{R_3,B_n}(1,q)^{\sharp}
\end{split}
\end{equation*}
Since $\sharp=\natural{\cdot}\flat$, we also have
\begin{equation*}
\begin{split}
\mathscr{H}_{R_3,B_n}(1,q)^{\natural}{\cap}\mathscr{H}_{R_3,B_n}(1,q)^{\sharp}
&=\mathscr{H}_{R_3,B_n}(1,q)^{\natural}{\cap}\mathscr{H}_{R_3,B_n}(1,q)^{\natural{\cdot}\flat}\\
&{\subseteq}\mathscr{H}_{R_3,B_n}(1,q)^{\natural}{\cap}\mathscr{H}_{R_3,B_n}(1,q)^{\flat}
\end{split}
\end{equation*}
Thus we obtain $\mathscr{H}_{R_3,B_n}(1,q)^{\natural}{\cap}\mathscr{H}_{R_3,B_n}(1,q)^{\flat}
=\mathscr{H}_{R_3,B_n}(1,q)^{\sharp}{\cap}\mathscr{H}_{R_3,B_n}(1,q)^{\natural}$. 
Other equations are also proved in the same fashion. 
\end{proof}
From the relations in Proposition 4.1, we can see that every monomial with occurrences of even 
(resp. odd) numbers of $a_1$ also has occurrences of even (resp. odd) numbers of $a_1$ 
in each term of any other expression of it in $\mathscr{H}_{R_3,B_n}(1,q)$, and that 
the same is holds for numbers of occurrences of $a_i$'s ($i>1$). 
Hence by the definition, $\mathscr{H}_{R_3,B_n}(1,q)^{\dagger}$ is generated by monomials with 
both occurrences of even numbers of $a_1$ and those of other $a_i$'s. 
We also have that 
\[
\mathscr{H}_{R_3,B_n}(1,q)^{\dagger}
=\bigoplus_{M\in\mathscr{E}_n^{\sharp}{\cap}\mathscr{E}_n^{\natural}}R_3M.
\]
If $n{\geq}4$, then $\mathscr{H}_{R_3,B_n}(1,q)^{\dagger}$ yields to be the 
subalgebra $\mathscr{H}_{R_3,D_n}(q)^{\natural}$ of 
$\mathscr{H}_{R_3,D_n}(q)$ consisting of fixed elements of $\natural$. 
\begin{Prop}
$\operatorname{rank}\mathscr{H}_{R_3,B_n}(1,q)^{\dagger}=2^{n-2}n!$. 
\end{Prop}
\begin{proof}
We show this by induction on $n$. 
For $n=2$, the assertion clearly holds.
Assume $n>2$. We define four subsets of $\mathscr{S}_n$ in 
$\mathscr{H}_{R_3,B_n}(1,q)$ as follows. 
\begin{equation*}
\begin{split}
X_n&=\{M\,|\,\text{$M$ has an even number of $b_1$ 
and an even number of $b_2,{\hdots},b_n$}\},\\
Y_n&=\{M\,|\,\text{$M$ has an odd number of $b_1$ 
and an even number of $b_2,{\hdots},b_n$}\},\\
Z_n&=\{M\,|\,\text{$M$ has an even number of $b_1$ 
and an odd number of $b_2,{\hdots},b_n$}\},\\
W_n&=\{M\,|\,\text{$M$ has an odd number of $b_1$ 
and an odd number of $b_2,{\hdots},b_n$}\}.
\end{split}
\end{equation*}
Write $M=U_1U_2{\cdots}U_n$. We consider four cases depending upon 
the subset to which $M'=U_1U_2{\cdots}U_{n-1}$ belongs. \\
case 1 : 
If $M'{\in}X_{n-1}$, then $U_n$ must have no $b_1$ and a 
even number of $b_2,{\hdots},b_n$. \\
case 2 : 
If $M'{\in}Y_{n-1}$, then $U_n$ must have one $b_1$ and a 
even number of $b_2,{\hdots},b_n$. \\
case 3 : 
If $M'{\in}Z_{n-1}$, then $U_n$ must have no $b_1$ and a 
odd number of $b_2,{\hdots},b_n$. \\
case 4 : 
If $M'{\in}W_{n-1}$, then $U_n$ must have one $b_1$ and a 
odd number of $b_2,{\hdots},b_n$. \\
If $n$ is even, exactly $n/2$ such $U_n$ exist in each case. 
Thus we obtain
\begin{equation*}
\begin{split}
|X_n|&=\dfrac{n}{2}|X_{n-1}|+\dfrac{n}{2}|Y_{n-1}|+\dfrac{n}{2}|Z_{n-1}|+\dfrac{n}{2}|W_{n-1}|\\
&=\dfrac{n}{2}|\mathscr{S}_{n-1}|=\dfrac{n}{2}2^{n-1}(n-1)!=2^{n-2}n!.
\end{split}
\end{equation*}
If $n$ is odd, then exactly $(n+1)/2$ such $U_n$ exist in case 1,2 and 
$(n-1)/2$ in the other cases. 
Hence, we have
\begin{equation*}
\begin{split}
|X_n|&=\dfrac{n+1}{2}|X_{n-1}|+\dfrac{n+1}{2}|Y_{n-1}|+\dfrac{n-1}{2}|Z_{n-1}|
+\dfrac{n-1}{2}|W_{n-1}|\\
&=\dfrac{n}{2}|\mathscr{S}_{n-1}|+\dfrac{1}{2}(|X_{n-1}|+|Y_{n-1}|-|Z_{n-1}|-|W_{n-1}|).
\end{split}
\end{equation*}
Observing $|X_{n-1}|+|Y_{n-1}|=|\mathscr{E}_{n-1}^{\flat}|=2^{n-2}(n-1)!$, we have 
$|Z_{n-1}|+|W_{n-1}|=2^{n-1}(n-1)!-2^{n-2}(n-1)!=2^{n-2}(n-1)!$. 
Thus $|X_n|=2^{-1}n|\mathscr{S}_{n-1}|=2^{n-2}n!$ as desired. 
\end{proof}
\begin{Theo}
$\mathscr{H}_{R_3,B_n}(1,q)^{\sharp}$, $\mathscr{H}_{R_3,B_n}(1,q)^{\flat}$ and 
$\mathscr{H}_{R_3,B_n}(1,q)^{\natural}$ are $\mathbb{Z}_2$-crossed products 
with $A_{\Bar{0}}=\mathscr{H}_{R_3,B_n}(1,q)^{\dagger}$. 
\end{Theo}
\begin{proof}
Let us define three submodules of $\mathscr{H}_{R_3,B_n}(1,q)$ as follows. 
\begin{equation*}
\begin{split}
\mathscr{H}_{R_3,B_n}(1,q)^{-\dagger,\sharp}&=\bigoplus_{M{\in}Y_n{\cup}W_n}R_3M,\\
\mathscr{H}_{R_3,B_n}(1,q)^{-\dagger,\flat}&=\bigoplus_{M{\in}Z_n{\cup}W_n}R_3M,\\
\mathscr{H}_{R_3,B_n}(1,q)^{-\dagger,\natural}&=\bigoplus_{M{\in}Y_n{\cup}Z_n}R_3M.
\end{split}
\end{equation*}
It is clear that there exist direct sum decompositions as follows. 
\begin{equation*}
\begin{split}
\mathscr{H}_{R_3,B_n}(1,q)^{\sharp}&=\mathscr{H}_{R_3,B_n}(1,q)^{\dagger}{\oplus}
\mathscr{H}_{R_3,B_n}(1,q)^{-\dagger,\sharp},\\
\mathscr{H}_{R_3,B_n}(1,q)^{\flat}&=\mathscr{H}_{R_3,B_n}(1,q)^{\dagger}{\oplus}
\mathscr{H}_{R_3,B_n}(1,q)^{-\dagger,\flat},\\
\mathscr{H}_{R_3,B_n}(1,q)^{\natural}&=\mathscr{H}_{R_3,B_n}(1,q)^{\dagger}{\oplus}
\mathscr{H}_{R_3,B_n}(1,q)^{-\dagger,\natural}.
\end{split}
\end{equation*}
One can readily check that these direct sum decompositions satisfy the 
condition of the $\mathbb{Z}_2$-crossed product. 
\end{proof}
\section{The branching rules for the Hecke algebras of Type $B$}
Let $K_0=\mathbb{Q}(u,q)$ be the quotient field of $R_0$ and $K_1=\mathbb{Q}(q)$ 
that of $R_1$. 
Let $\lambda=(\lambda^{(1)},\lambda^{(2)})$ be a $2$-tuple of Young diagrams of total 
size $n$. Let $T=(T^{(1)},T^{(2)})$ be a tableau of shape $\lambda$. 
We mean by a tableau $T=(T^{(1)},T^{(2)})$ of shape $\lambda=(\lambda^{(1)},\lambda^{(2)})$ 
the shape of $T^{(i)}$ is $\lambda^{(i)}(i=1,2)$, and each of the symbols $1,2,\hdots ,n$ 
appears in $T$ exactly once. 
For $\lambda=(\lambda^{(1)},\lambda^{(2)})$(resp. $T=(T^{(1)},T^{(2)})$), let 
${\lambda^{(i)}{'}}$(resp. ${T^{(i)}{'}}$) be its transpose for $i=1,2$, and 
$\lambda'=(\lambda^{(1)}{'},\lambda^{(2)}{'})$(resp. $T'=(T^{(1)}{'},T^{(2)}{'})$). 
Moreover, Let 
$\lambda^*=(\lambda^{(2)},\lambda^{(1)})$, $T^*=(T^{(2)},T^{(1)})$ and 
$\lambda'^*=(\lambda^{(2)}{'},\lambda^{(1)}{'})$, $T'^*=(T^{(2)}{'},T^{(1)}{'})$. 
A tableau of shape $\lambda=(\lambda^{(1)},\lambda^{(2)})$ 
is said to be \it standard \rm if the numbers $1,2,\hdots,n$ increase along the rows and 
columns of each Young diagram $\lambda^{(1)},\lambda^{(2)}$ of $\lambda$. 
We denote by $T^{(\alpha)}(i,j)(\alpha=1,2)$ the number at $i$-th 
row and $j$-th column of $T^{(\alpha)}$. 
Let $\operatorname{STab}(\lambda)$ be the set of all of standard tableaux of shape $\lambda$. 
Let $k$ be the number at $i$-th row and $j$-th column of $T^{(1)}$ or $T^{(2)}$ of 
$T{\in}\operatorname{STab}(\lambda)$, and $l$ at $i'$-th row and $j'$-th column of them. 
Then the integer $d_{T,k,l} = (j'-i')-(j-i)$ is said to be 
\it the axial distance from $k$ to $l$ in $T$. \rm
For each $T{\in}\operatorname{STab}(\lambda)$ we take a symbol $v_T$, and 
define the free $K_0$-module 
\[ V_\lambda  = \bigoplus_{T\in\operatorname{STab}(\lambda)}K_0v_{T}. \]
For each $\lambda$, we can give an irreducible representation of $\mathscr{H}_{K_0,B_n}(u,q)$ 
in the following manner. 
Let $y$ be any indeterminate and $k$ any integer. 
Let $M(k,y)$ be the $2{\times}2$ matrix defined by
\[
M(k,y)=\dfrac{1}{1-q^ky}
	\begin{bmatrix}
		q-1 & 1-q^{k+1}y \\
		q(1-q^{k-1}y) & -q^ky(q-1)
	\end{bmatrix}
\]
For a standard tableau $T=(T^{(1)},T^{(2)})$ of shape 
$\lambda=(\lambda^{(1)},\lambda^{(2)})$ of total size $n$, we define $\rho_T$ to be the map 
from $\{1,2,\hdots,n\}$ to $\{1,2\}$ such that the number $i$ 
occurs in the $\rho_T(i)$-th Young diagram $T^{(\rho_T(i))}$ of $T$. 
We set $u_1=u$ and $u_2=-1$. \par
We shall give an action of the generators of $\mathscr{H}_{K_0,B_n}(u,q)$ on $V_\lambda$ 
as follows. 
\renewcommand{\theenumi}{\arabic{enumi}}
\renewcommand{\labelenumi}{(a\theenumi)} 
\begin{enumerate}
\item
$a_1v_{T}=u_{\rho_T(1)}v_{T}$.
\item
If $i>1$, $a_i$ acts on $V_\lambda$ in three ways, depending upon the 
position that $i-1$ and $i$ occupy in $T=(T^{(1)},T^{(2)})$ as follows.
\begin{enumerate}
\item
If $i-1$ and $i$ appear in the same row of the same diagram of $(T^{(1)},T^{(2)})$, then
$a_iv_{T}=qv_{T}$.
\item
If $i-1$ and $i$ appear in the same column of the same diagram of $(T^{(1)},T^{(2)})$, then
$a_iv_{T}=-v_{T}$.
\item
Elsewhere, $a_i$ acts on the subspace $K_0v_{T}\oplus K_0v_{s_{i-1}T}$ 
of $V_\lambda$ as follows, 
where $s_{i-1}=(i-1,i)$ is the transposition of $i-1$ and $i$, and 
$s_{i-1}T$ is the standard tableau obtained from $T$ 
by transposing $i-1$ and $i$. 
\[
a_i\langle v_{T},v_{s_{i-1}T}\rangle = \langle v_{T},v_{s_{i-1}T}\rangle 
	M\Big{(}d_{T,i,i-1},\dfrac{u_{\rho_T(i-1)}}{u_{\rho_T(i)}}\Big{)}.
\]
\end{enumerate}
\end{enumerate}
\renewcommand{\theenumi}{\arabic{enumi}}
\renewcommand{\labelenumi}{(\theenumi)}
\begin{Theo}[\cite{Hoefsmit},Theorem 2.2.11]
For each $2$-tuple of Young diagrams of total size $n$, the above action of 
$\mathscr{H}_{K_0,B_n}(u,q)$ gives an (absolutely) irreducible 
representation of $\mathscr{H}_{K_0,B_n}(u,q)$. If $\lambda \neq \mu$ as $2$-tuples of 
Young diagrams, then $V_\lambda$ and $V_\mu$ are mutually inequivalent irreducible 
representations of $\mathscr{H}_{K_0,B_n}(u,q)$. These constitute 
a complete set of representatives of irreducible representations of 
$\mathscr{H}_{K_0,B_n}(u,q)$
\end{Theo}
These representations are said to be the \it seminormal form representations \rm 
of $\mathscr{H}_{K_0,B_n}(u,q)$. 
We denote by $(\pi_\lambda,V_\lambda)$ the seminormal form representation of
$\mathscr{H}_{K_0,B_n}(u,q)$ corresponding to $\lambda$. 
By an easy calculation, we obtain that $\pi_{\lambda}(b_i)$ is as follows.
\renewcommand{\theenumi}{\arabic{enumi}}
\renewcommand{\labelenumi}{(b\theenumi)} 
\begin{enumerate}
\item
$\pi_\lambda(b_1)v_{T}=(-1)^{\rho_T(1)-1}v_{T}$.
\item
If $i>1$, $\pi_\lambda(b_i)$ are given in three ways, depending upon the 
position that $i-1$ and $i$ occupy in $T=(T^{(1)},T^{(2)})$ as follows.
\begin{enumerate}
\item
If $i-1$ and $i$ appear in the same row of the same diagram of $(T^{(1)},T^{(2)})$, then
$\pi_\lambda(b_i)v_{T}=v_{T}$.
\item
If $i-1$ and $i$ appear in the same column of the same diagram of $(T^{(1)},T^{(2)})$, then
$\pi_\lambda(b_i)v_{T}=-v_{T}$.
\item
Elsewhere, $\pi_\lambda(b_i)$ acts on the subspace 
$K_0v_{T}\oplus K_0v_{s_{i-1}T}$ of $V_\lambda$ as follows. 
\[
\pi_{\lambda}(b_i){\langle}v_T,v_{s_{i-1}T}{\rangle}
={\langle}v_T,v_{s_{i-1}T}{\rangle} 
M'\Big{(}d_{T,i,i-1},\dfrac{u_{\rho_T(i-1)}}{u_{\rho_T(i)}}\Big{)}
\]
where $M'(k,y)$ is a $2\times 2$ matrix defined by
\[
M'(k,y)=\dfrac{1}{(q+1)(1-q^ky)}
	\begin{bmatrix}
		(q-1)(1+q^ky) & 2(1-q^{k+1}y) \\
		2q(1-q^{k-1}y) & -(q-1)(1+q^ky)
	\end{bmatrix}
\]
\end{enumerate}
\end{enumerate}
\renewcommand{\theenumi}{\arabic{enumi}}
\renewcommand{\labelenumi}{(\theenumi)}
We denote by $(\pi_{\lambda,\natural},V_{\lambda,\natural})$ the restriction of 
$(\pi_\lambda,V_\lambda)$ to $\mathscr{H}_{K_0,B_n}(u,q)^{\natural}$. 
We obtained the branching rule between $\mathscr{H}_{K_0,B_n}(u,q)$ and 
$\mathscr{H}_{K_0,B_n}(u,q)^{\natural}$ in \cite{Mitsuhashi2} as follows. 
\begin{Theo}[\cite{Mitsuhashi2} Corollary 4.5, Proposition 4.7, Theorem 4.8]
\begin{enumerate}
\rm \item \it
$\pi_{\lambda'^*,\natural}{\cong}\pi_{\lambda,\natural}$.
\rm \item \it
If $\lambda'^*={\lambda}$, then $\pi_{\lambda,\natural}$ decomposes into 
inequivalent subrepresentations 
$\pi_{\lambda,\natural}^+$ and $\pi_{\lambda,\natural}^-$ over $\Bar{K}_0$, 
the algebraic closure of $K_0$. Furthermore, $\deg\pi_{\lambda,\natural}^+
=\deg\pi_{\lambda,\natural}^-=\deg\pi_{\lambda,\natural}/2$.
\rm \item \it
Let $\{\lambda_i,\lambda_i'^*,\mu_j\}_{i,j}$ 
be the set of $2$-tuples of Young diagrams 
such that 
\begin{equation*}
\lambda_i'^*{\neq}\lambda_i,{\quad}\mu_j'^*=\mu_j
\end{equation*}
Then 
\[
\operatorname{Irr}(\mathscr{H}_{\bar{K}_0,B_n}(q)^{\natural})=
\{\pi_{\lambda_i,\natural},\pi_{\mu_j,\natural}^+,\pi_{\mu_j,\natural}^-,\}_{i,j}
\]
is a basic set of irreducible representations of $\mathscr{H}_{\bar{K}_0,B_n}(u,q)^{\natural}$. 
$\mathscr{H}_{\bar{K}_0,B_n}(u,q)^{\natural}$ is semisimple. 
\end{enumerate}
\end{Theo}
The proof is in \cite{Mitsuhashi2}, but there are some errors in it. 
So we shall prove again. 
From (b1),(b2), we obtain that $\pi_{\lambda'^*}(b_i)$ is as follows by 
a direct computation. 
\renewcommand{\theenumi}{\arabic{enumi}}
\renewcommand{\labelenumi}{(c\theenumi)}
\begin{enumerate}
\item
$\pi_{\lambda'^*}(b_1)v_{T'^*}=(-1)^{\rho_T(1)}v_{T'^*}$.
\item
If $i>1$, $\pi_{\lambda'^*}(b_i)$ are given in three ways, depending upon the 
position that $i-1$ and $i$ occupy in $T=(T^{(1)},T^{(2)})$ as follows.
\begin{enumerate}
\item
If $i-1$ and $i$ appear in the same row of the same diagram of $T$, then
$\pi_{\lambda'^*}(b_i)v_{T'^*}=-v_{T'^*}$.
\item
If $i-1$ and $i$ appear in the same column of the same diagram of $T$, then
$\pi_{\lambda'^*}(b_i)v_{T'^*}=v_{T'^*}$.
\item
Elsewhere, $\pi_{\lambda'^*}(b_i)$ reduces to the endomorphism of the subspace 
$K_0v_{T'^*}\oplus K_0v_{s_{i-1}{T'^*}}$ of $V_{\lambda'^*}$ as follows. 
\begin{equation*}
\begin{split}
\pi_{\lambda'^*}(b_i){\langle}v_{T'^*},v_{s_{i-1}{T'^*}}{\rangle}
&={\langle}v_{T'^*},v_{s_{i-1}{T'^*}}{\rangle} 
M'\Big{(}-d_{T,i,i-1},\dfrac{u_{\rho_{T}(i)}}{u_{\rho_{T}(i-1)}}\Big{)}\\
&={\langle}v_{T'^*},v_{s_{i-1}{T'^*}}{\rangle} 
JM'\Big{(}d_{T,i,i-1},\dfrac{u_{\rho_{T}(i-1)}}{u_{\rho_{T}(i)}}\Big{)}J
\end{split}
\end{equation*}
where $J=\begin{bmatrix}0&1\\1&0\end{bmatrix}$. 
\end{enumerate}
\end{enumerate}
\renewcommand{\theenumi}{\arabic{enumi}}
\renewcommand{\labelenumi}{(\theenumi)}
In \cite{Mitsuhashi2}, p.243-244, the intertwining operator between 
$\pi_{\lambda'^*,\natural}$ and $\pi_{\lambda,\natural}$ is given. 
But it is incorrect since the denominator of $\alpha_T(i,j)$ equals $0$ 
if $i$ and $j$ belong to the same Young diagram and satisfy $d_{T,i,j}=0$. 
So we correct the intertwining operator as follows. 
For each standard tableau $T$ of shape $\lambda$ of total size $n$, 
we define the map $\psi_T$ from 
$I=\{(i,j)\in \mathbb{N}\times\mathbb{N}\mid 1\leq i{\neq}j \leq n \}$ to $K_0$ 
divided into two cases as follows. \\
case 1 : $i$ and $j$ appear in the same Young diagram.
\begin{equation*}
\psi_T(i,j)=
\begin{cases}
	1 \qquad \text{if $i$ and $j$ appear in the same row or the same column,}\\
	\hspace{1cm}\text{or $i$ and $j$ satisfy $d_{T,i,j}=0$,} \\
	\dfrac{q(1-q^{d_{T,i,j}-1}u_{\rho_T(j)}u_{\rho_T(i)}^{-1})}
		{(q+1)(1-q^{d_{T,i,j}}u_{\rho_T(j)}u_{\rho_T(i)}^{-1})} \times 
		\operatorname{sgn}(d_{T,i,j}) \qquad \text{otherwise.} 
\end{cases}
\end{equation*}
case 2 : $i$ and $j$ appear in the different Young diagrams.
\begin{equation*}
\psi_T(i,j)=
	\dfrac{q(1-q^{d_{T,i,j}-1}u_{\rho_T(j)}u_{\rho_T(i)}^{-1})}
		{(q+1)(1-q^{d_{T,i,j}}u_{\rho_T(j)}u_{\rho_T(i)}^{-1})} \times 
		\Big{(}\rho_T(i)-\rho_T(j)\Big{)}
\end{equation*}
Let $\psi:\operatorname{STab}(\lambda)\longrightarrow{K_0}$ be the map defined by 
\begin{equation}
\psi(T)=\prod_{1{\leq}j<i{\leq}n}\psi_T(i,j),
\end{equation}
and $\varPsi_{\natural}$ the $K_0$-homomorphism defined by
\begin{equation*}
\begin{split}
\varPsi_{\natural} : \quad &V_{\lambda} \longrightarrow V_{\lambda'^*} \\
	&v_T \longmapsto \psi(T)v_{T'^*} \\
\end{split}
\end{equation*}
Then, Proposition 4.4 in \cite{Mitsuhashi2} holds correctly as follows. 
\begin{Prop}
\begin{equation*}
\begin{split}
\pi_{\lambda'^*}(b_i)\varPsi_{\natural}(v_{T}) 
	&=-\sum_{S\in \operatorname{STab}(\lambda)}
	\Big{(}\pi_\lambda(b_i)\Big{)}_{T,S}\varPsi_{\natural}(v_{S}) \qquad 
	\text{for $i=1,2,\hdots,n$}  
\end{split}
\end{equation*}
\end{Prop}
\begin{proof}
In the cases (c1) and (c2)-(a)(b), the equation holds obviously. 
Hence we consider only the case (c2)-(c). 
By a direct calculation, we have 
\begin{equation*}
\begin{split}
&\dfrac{(1-q^{d_{T,i,i-1}+1}u_{\rho_T(i-1)}u_{\rho_T(i)}^{-1})}
	{(q+1)(1-q^{d_{T,i,i-1}}u_{\rho_T(i-1)}u_{\rho_T(i)}^{-1})}\\
&=\dfrac{q(1-q^{d_{s_{i-1}T,i,i-1}-1}u_{\rho_{s_{i-1}T}(i-1)}u_{\rho_{s_{i-1}T}(i)}^{-1})}
	{(q+1)(1-q^{d_{s_{i-1}T,i,i-1}}u_{\rho_{s_{i-1}T}(i-1)}u_{\rho_{s_{i-1}T}(i)}^{-1})}
\end{split}
\end{equation*}
Hence by the definition of $\psi(T)$, one can check the following. 
\begin{equation*}
\begin{split}
&\dfrac{(1-q^{d_{T,i,i-1}+1}u_{\rho_T(i-1)}u_{\rho_T(i)}^{-1})}
	{(q+1)(1-q^{d_{T,i,i-1}}u_{\rho_T(i-1)}u_{\rho_T(i)}^{-1})}\psi(T)\\
&=-\dfrac{q(1-q^{d_{T,i,i-1}-1}u_{\rho_T(i-1)}u_{\rho_T(i)}^{-1})}
	{(q+1)(1-q^{d_{T,i,i-1}}u_{\rho_T(i-1)}u_{\rho_T(i)}^{-1})}\psi(s_{i-1}T)
\end{split}
\end{equation*}
Therefore, we obtain the following equation.
\begin{equation*}
\begin{split}
\pi_{\lambda'^*}(b_i)\psi(T)v_{T'^*}
	&=\dfrac{(q-1)(1+q^{d_{T'^*,i,i-1}}u_{\rho_{T'^*}(i-1)}u_{\rho_{T'^*}(i)}^{-1})}
		{(q+1)(1-q^{d_{T'^*,i,i-1}}u_{\rho_{T'^*}(i-1)}u_{\rho_{T'^*}(i)}^{-1})}
		\psi(T)v_{T'^*}\\
	&+\dfrac{2q(1-q^{d_{T'^*,i,i-1}-1}u_{\rho_{T'^*}(i-1)}u_{\rho_{T'^*}(i)}^{-1})}
		{(q+1)(1-q^{d_{T'^*,i,i-1}}u_{\rho_{T'^*}(i-1)}u_{\rho_{T'^*}(i)}^{-1})}
		\psi(T)v_{s_{i-1}T'^*}\\
	&=-\dfrac{(q-1)(1+q^{d_{T,i,i-1}}u_{\rho_T(i-1)}u_{\rho_T(i)}^{-1})}
		{(q+1)(1-q^{d_{T,i,i-1}}u_{\rho_T(i-1)}u_{\rho_T(i)}^{-1})}
		\psi(T)v_{T'^*}\\
	&-\dfrac{2q(1-q^{d_{T,i,i-1}-1}u_{\rho_T(i-1)}u_{\rho_T(i)}^{-1})}
		{(q+1)(1-q^{d_{T,i,i-1}}u_{\rho_T(i-1)}u_{\rho_T(i)}^{-1})}
		\psi(s_{i-1}T)v_{s_{i-1}T'^*}
\end{split}
\end{equation*}
Thus we have proved the proposition. 
\end{proof}
\begin{Cor}
$\pi_{\lambda,\natural}{\cong}\pi_{\lambda'^*,\natural}$. 
The $K_0$-homomorphism 
$\varPsi_{\natural}:V_{\lambda,\natural}{\longrightarrow}V_{\lambda'^*,\natural}$ 
is the intertwining operator between $\pi_{\lambda,\natural}$ and 
$\pi_{\lambda'^*,\natural}$. 
\end{Cor}
\begin{proof}
$\mathscr{H}_{K_0,B_n}(u,q)^{\natural}$ is generated by 
all the monomials of even numbers of products of $b_1,b_2,{\hdots},b_n$. 
Thus the representation matrix of 
$\pi_{\lambda,\natural}$ with respect to the basis 
$\{v_T{\mid}T{\in}\operatorname{STab}(\lambda)\}$ and that of 
$\pi_{\lambda'^*,\natural}$ with respect to the basis 
$\{\varPsi_{\natural}(v_{T}){\mid}T{\in}\operatorname{STab}(\lambda)\}$ coincide.
\end{proof}
Therefore we obtain (1) of Theorem 6.2. 
\begin{Cor}
Let $b{\in}\mathscr{H}_{K_0,B_n}(u,q)^{\natural}$. Then
\[ \Big{(}\pi_{\lambda'^*}(b)\Big{)}_{T'^*,S'^*}\psi(T)
	= \Big{(}\pi_{\lambda}(b)\Big{)}_{T,S}\psi(S)
\]
\end{Cor}
\begin{proof}
Let 
$b=b_{j_1}b_{j_2}\cdots b_{j_k}{\in}\mathscr{H}_{K_0,B_n}(u,q)^{\natural}$. 
Let $T_1=T,T_{k}=S$. Then we may write 
\[ \pi_{\lambda'^*}(b)\psi(T_1)v_{T_1'^*} 
	= \sum_{T_{k}{\in}\operatorname{STab}(\lambda)}
	\Big{(}\pi_{\lambda'^*}(b)\Big{)}_{T_1'^*,T_k'^*}
	\psi(T_1)v_{T_k'^*}.
\]
Since $k$ is even, we obtain the following by using Proposition 6.3 repeatedly. 
\begin{equation*}
\begin{split}
&\pi_{\lambda'^*}(b)\psi(T_1)v_{T_1'^*} 
	=\pi_{\lambda'^*}(b_{j_1})\pi_{\lambda'^*}(b_{j_2})\cdots 
	\pi_{\lambda'^*}(b_{j_k})	\psi(T_1)v_{T_1'^*} \\
	&= \sum_{T_2,T_3,\hdots,T_k{\in}\operatorname{STab}(\lambda)}
	\Big{(}\pi_{\lambda}(b_{j_k})\Big{)}_{T_1,T_2}
	\Big{(}\pi_{\lambda}(b_{j_{k-1}})\Big{)}_{T_2,T_3}\cdots 
	\Big{(}\pi_{\lambda}(b_{j_1})\Big{)}_{T_{k-1},T_k}
	\psi(T_k)v_{T_k'^*} \\
	&= \sum_{T_k{\in}\operatorname{STab}(\lambda)}
	\Big{(}\pi_{\lambda}(b)\Big{)}_{T_1,T_k}\psi(T_k)v_{T_k'^*}.
\end{split}
\end{equation*}
Comparing the coefficients of $v_{T_k'^*}$ of the both right-hand sides, 
we obtain the assertion. 
\end{proof}
For $\lambda$ such that $\lambda'^*=\lambda$, 
we define $\operatorname{STab(\lambda)}^{\natural}$ and 
$\operatorname{STab(\lambda)}^{-\natural}$ to be the sets of all the standard tableaux 
of shape $\lambda$ such that $1$ belongs to $T^{(1)}$ and $T^{(2)}$ respectively. 
Then, $\operatorname{STab(\lambda)}$ is the disjoint union of 
$\operatorname{STab(\lambda)}^{\natural}$ and $\operatorname{STab(\lambda)}^{-\natural}$. 
There exists an involutive bijection from $\operatorname{STab(\lambda)}^{\natural}$ to 
$\operatorname{STab(\lambda)}^{-\natural}$ which is defined by $T{\mapsto}T'^*$. 
The submodules $V_{\lambda,\natural}^+,V_{\lambda,\natural}^-$ of $V_{\lambda,\natural}$ 
have also been given in p.246 of \cite{Mitsuhashi2} as follows. 
\begin{equation*}
\begin{split}
V_{\lambda,\natural}^+ &= \bigoplus_{T{\in}\operatorname{STab}(\lambda)^{\natural}}
\Bar{K}_1(\sqrt{\psi(T'^*)}v_{T}+\sqrt{\psi(T)}v_{T'^*}), \\
V_{\lambda,\natural}^- &= \bigoplus_{T{\in}\operatorname{STab}(\lambda)^{\natural}}
\Bar{K}_1(\sqrt{\psi(T'^*)}v_{T}-\sqrt{\psi(T)}v_{T'^*}),
\end{split}
\end{equation*}
We notice that $\sqrt{\psi(T)}$ has two branches. 
We take a suitable branch in each computation of square roots for consistency. 
As in the proof of Proposition 4.7 of \cite{Mitsuhashi2}, we also have 
\begin{equation*}
\begin{split}
&\pi_{\lambda,\natural}(b)\Big{(}\sqrt{\psi(T'^*)}v_{T}{\pm}\sqrt{\psi(T)}v_{T'^*}\Big{)} 
	= \sum_{S{\in}\operatorname{STab}(\lambda)^{\natural}}\Big{\{}
		\sqrt{\Big{(}\pi_{\lambda,\natural}(b)\Big{)}_{T,S}
		\Big{(}\pi_{\lambda,\natural}(b)\Big{)}_{T'^*,S'^*}} \\
		&{\quad}{\pm}\sqrt{\Big{(}\pi_{\lambda,\natural}(b)\Big{)}_{T,S'^*} 
		\Big{(}\pi_{\lambda,\natural}(b)\Big{)}_{T'^*,S}}\Big{\}} 
		\times\Big{(}\sqrt{\psi(S'^*)}v_{S}
		+\sqrt{\psi(S)}v_{S'^*}\Big{)} \\
\end{split}
\end{equation*}
for $T{\in}\operatorname{STab}(\lambda)^{\natural}$. 
\begin{proof}[proof of Theorem 6.2]
(1) is Corollary 6.4 itself. 
If $\lambda'^*=\lambda$, then one can deduce directly from Theorem 2.7 that 
both $\pi_{\lambda,\natural}^+$ and $\pi_{\lambda,\natural}^-$ are irreducible and 
mutually inequivalent. Thus (2) holds. 
Theorem 2.7 also shows immediately that 
$\pi_{\lambda_i,\natural},\pi_{\mu_j,\natural}^+,\pi_{\mu_j,\natural}^-$ are 
mutually inequivalent. 
Semisimplicity of $\mathscr{H}_{\bar{K}_0,B_n}(u,q)^{\natural}$ has been 
proved in Theorem 4.8 of \cite{Mitsuhashi2}. 
\end{proof}
On the other hand, The branching rule between $\mathscr{H}_{K_1,B_n}(1,q)$ and 
$\mathscr{H}_{K_1,B_n}(1,q)^{\sharp}$ has given in \cite{Hoefsmit} for $n{\geq}4$. 
We also refer \cite{G-P} for the detail of the irreducible representations of 
$\mathscr{H}_{K_1,D_n}(q)$. 
We shall give a proof of the branching rule for $n{\geq}2$ using the theory of 
crossed products. Now we set $u=1$. 
We denote by $(\pi_{\lambda,\sharp},V_{\lambda,\sharp})$ the restriction of 
$(\pi_\lambda,V_\lambda)$ to $\mathscr{H}_{K_1,B_n}(1,q)^{\sharp}$. 
From (b1),(b2), we obtain that $\pi_{\lambda^*}(b_i)$ is as follows. 
\renewcommand{\theenumi}{\arabic{enumi}}
\renewcommand{\labelenumi}{(d\theenumi)}
\begin{enumerate}
\item
$\pi_{\lambda^*}(b_1)v_{T^*}=(-1)^{\rho_T(1)}v_{T^*}$.
\item
If $i>1$, $\pi_{\lambda^*}(b_i)$ are given in three ways, depending upon the 
position that $i-1$ and $i$ occupy in $T=(T^{(1)},T^{(2)})$ as follows.
\begin{enumerate}
\item
If $i-1$ and $i$ appear in the same row of the same diagram of $T$, then
$\pi_{\lambda^*}(b_i)v_{T^*}=v_{T^*}$.
\item
If $i-1$ and $i$ appear in the same column of the same diagram of $T$, then
$\pi_{\lambda^*}(b_i)v_{T^*}=-v_{T^*}$.
\item
Elsewhere, $\pi_{\lambda^*}(b_i)$ acts on the subspace 
$K_0v_{T^*}\oplus K_0v_{s_{i-1}{T^*}}$ of $V_{\lambda^*}$ as follows. 
\[
\pi_{\lambda^*}(b_i){\langle}v_{T^*},v_{s_{i-1}{T^*}}{\rangle}
={\langle}v_{T^*},v_{s_{i-1}{T^*}}{\rangle} 
M'\Big{(}d_{T,i,i-1},\dfrac{u_{\rho_{T}(i-1)}}{u_{\rho_{T}(i)}}\Big{)}
\]
\end{enumerate}
\end{enumerate}
\renewcommand{\theenumi}{\arabic{enumi}}
\renewcommand{\labelenumi}{(\theenumi)}
Indeed, (d1) and (d2)-(a)(b) are obvious, and (d2)-(c) follows from 
$d_{T^*,i.i-1}=d_{T,i.i-1},
\,u_{\rho_{T^*}(i-1)}/u_{\rho_{T^*}(i)}=u_{\rho_{T}(i-1)}/u_{\rho_{T}(i)}$. 
Therefore, (d1),(d2) give the following. 
\begin{Prop}[\cite{G-P} 10.4]
$\pi_{\lambda,\sharp}{\cong}\pi_{\lambda^*,\sharp}$. Especially, two matrices 
correspond to $\pi_{\lambda^*,\sharp}(h)$ and $\pi_{\lambda,\sharp}(h)$ 
($h{\in}\mathscr{H}_{K_1,B_n}(1,q)^{\sharp}$) coincide. 
\end{Prop}
\begin{proof}
Since $\mathscr{H}_{K_1,B_n}(1,q)^{\sharp}$ is generated by 
all the monomials with occurrences of even numbers of $b_1$, 
all the matrices $\pi_{\lambda^*,\sharp}(h)$ ($h{\in}\mathscr{H}_{K_1,B_n}(1,q)^{\sharp}$) 
coincide with $\pi_{\lambda,\sharp}(h)$. 
\end{proof}
If $\lambda^*=\lambda$, we shall show that $V_{\lambda,\sharp}$ decomposes into 
two nonzero submodules. 
Let $\omega$ be the endomorphism of $V_{\lambda}$ defined by 
\[
\omega(v_T)=v_{T^*}
\]
Then, obviously $\omega^2=1$. Furthermore, $\omega$ satisfies the following property. 
\begin{Prop}
$\omega\pi_{\lambda}(h)=\pi_{\lambda}(h)\omega$ for all 
$h{\in}\mathscr{H}_{K_1,B_n}(1,q)^{\sharp}$. 
\end{Prop}
\begin{proof}
If $h=b_1$, then
\begin{equation*}
\begin{split}
\omega\big{(}\pi_{\lambda}(b_1)v_T\big{)}
&=(-1)^{\rho_{T}(1)-1}{\omega}(v_T)=(-1)^{\rho_{T}(1)-1}v_{T^*}
=-(-1)^{\rho_{T}(1)}v_{T^*}\\
\pi_{\lambda}(b_1){\omega}(v_T)&=\pi_{\lambda}(b_1)v_{T^*}=(-1)^{\rho_{T}(1)}v_{T^*}
\end{split}
\end{equation*}
Hence $\omega\big{(}\pi_{\lambda}(b_1)v_T\big{)}=-\pi_{\lambda}(b_1){\omega}(v_T)$. \\
If $h=b_i$ ($i>1$), then by (b2) and (d2) we have
\[
\omega\big{(}\pi_{\lambda}(b_i)v_T\big{)}=\pi_{\lambda}(b_i){\omega}(v_T).
\]
Since $\mathscr{H}_{K_1,B_n}(1,q)^{\sharp}$ is generated by 
all the monomials with occurrences of even numbers of $b_1$, 
$\omega\pi_{\lambda}(h)=\pi_{\lambda}(h)\omega$ holds. 
\end{proof}
Let $\lambda^*=\lambda$. 
We define $\operatorname{STab(\lambda)}^{\sharp}$ and $\operatorname{STab(\lambda)}^{-\sharp}$ 
to be the sets of all the standard tableaux 
$T$ of shape $\lambda$ such that $1$ belongs to $T^{(1)}$ and $T^{(2)}$ respectively. 
Clearly $\operatorname{STab(\lambda)}$ is the disjoint union of 
$\operatorname{STab(\lambda)}^{\sharp}$ and $\operatorname{STab(\lambda)}^{-\sharp}$. 
There exists an involutive bijection from $\operatorname{STab(\lambda)}^{\sharp}$ to 
$\operatorname{STab(\lambda)}^{-\sharp}$ which is defined by $T{\mapsto}T^*$. 
We define two subspaces $V_{\lambda,\sharp}^+,V_{\lambda,\sharp}^+$ as follows. 
\begin{equation*}
\begin{split}
V_{\lambda,\sharp}^+&=\bigoplus_{T{\in}\operatorname{STab}(\lambda)^{\sharp}}K_1(v_T+v_{T^*})\\
V_{\lambda,\sharp}^-&=\bigoplus_{T{\in}\operatorname{STab}(\lambda)^{\sharp}}K_1(v_T-v_{T^*})
\end{split}
\end{equation*}
Then ${\dim}V_{\lambda,\sharp}^+={\dim}V_{\lambda,\sharp}^-={\dim}V_{\lambda}/2$ and 
$V_{\lambda}=V_{\lambda,\sharp}^+{\oplus}V_{\lambda,\sharp}^-$. Moreover, 
$V_{\lambda,\sharp}^+$ and $V_{\lambda,\sharp}^-$ are the eigenspaces corresponding 
to the eigenvalues $1$ and $-1$ of $\omega$ respectively. 
Therefore, $V_{\lambda,\sharp}^+$ and $V_{\lambda,\sharp}^-$ are 
$\mathscr{H}_{K_1,B_n}(1,q)^{\sharp}$-submodules. We denote by $(\pi_{\lambda,\sharp}^+,V_{\lambda,\sharp}^+)$ and 
$(\pi_{\lambda,\sharp}^-,V_{\lambda,\sharp}^-)$ the representations corresponding to the 
submodules $V_{\lambda,\sharp}^+,V_{\lambda,\sharp}^-$ respectively. 
\begin{Theo}
\begin{enumerate}
\rm \item \it
$\pi_{\lambda^*,\sharp}{\cong}\pi_{\lambda,\sharp}$.
\rm \item \it
If $\lambda^*={\lambda}$, then $\pi_{\lambda,\sharp}$ decomposes into 
inequivalent subrepresentations 
$\pi_{\lambda,\sharp}^+$ and $\pi_{\lambda,\sharp}^-$ over $\Bar{K}_1$, 
the algebraic closure of $K_1$. Furthermore, $\deg\pi_{\lambda,\sharp}^+
=\deg\pi_{\lambda,\sharp}^-=\deg\pi_{\lambda,\sharp}/2$.
\rm \item \it
Let $\{\lambda_i,\lambda_i^*,\mu_j\}_{i,j}$ 
be the set of $2$-tuples of Young diagrams 
such that 
\begin{equation*}
\lambda_i^*{\neq}\lambda_i,{\quad}\mu_j^*=\mu_j
\end{equation*}
Then 
\[
\operatorname{Irr}(\mathscr{H}_{\bar{K}_1,B_n}(q)^{\sharp})=
\{\pi_{\lambda_i,\sharp},\pi_{\mu_j,\sharp}^+,\pi_{\mu_j,\sharp}^-,\}_{i,j}
\]
is a basic set of irreducible representations of $\mathscr{H}_{\bar{K}_1,B_n}(q)^{\sharp}$. 
$\mathscr{H}_{\bar{K}_1,B_n}(q)^{\sharp}$ is semisimple. 
\end{enumerate}
\end{Theo}
\begin{proof}
(1) is Proposition 6.6 itself. 
The submodules $V_{\lambda,\sharp}^+$ and $V_{\lambda,\sharp}^-$ 
correspond to the case (2) of Theorem 2.7. Hence they are simple and mutually 
non-isomorphic. Thus (2) holds. 
The proof of (3) is in the same manner as in the proof of Theorem 6.2, 
just replacing $\natural$ and $\bar{K}_0$ with $\sharp$ and $\bar{K}_1$ respectively. 
So we omit the detail. 
\end{proof}
Next we shall give the branching rule between 
$\mathscr{H}_{K_1,B_n}(1,q)$ and $\mathscr{H}_{K_1,B_n}(1,q)^{\flat}$. 
From (b1),(b2), we obtain that $\pi_{\lambda'}(b_i)$ is as follows. 
\renewcommand{\theenumi}{\arabic{enumi}}
\renewcommand{\labelenumi}{(e\theenumi)}
\begin{enumerate}
\item
$\pi_{\lambda'}(b_1)v_{T'}=(-1)^{\rho_T(1)-1}v_{T'}$.
\item
If $i>1$, $\pi_{\lambda'}(b_i)$ are given in three ways, depending upon the 
position that $i-1$ and $i$ occupy in $T=(T^{(1)},T^{(2)})$ as follows.
\begin{enumerate}
\item
If $i-1$ and $i$ appear in the same row of the same diagram of $T$, then
$\pi_{\lambda'}(b_i)v_{T'}=-v_{T'}$.
\item
If $i-1$ and $i$ appear in the same column of the same diagram of $T$, then
$\pi_{\lambda'}(b_i)v_{T'}=v_{T'}$.
\item
Elsewhere, $\pi_{\lambda'}(b_i)$ acts on the subspace 
$K_0v_{T'}\oplus K_1v_{s_{i-1}{T'}}$ of $V_{\lambda'}$ as follows. 
\[
\pi_{\lambda'}(b_i){\langle}v_{T'},v_{s_{i-1}{T'}}{\rangle}
={\langle}v_{T'},v_{s_{i-1}{T'}}{\rangle} 
M'\Big{(}-d_{T,i,i-1},\dfrac{u_{\rho_{T}(i-1)}}{u_{\rho_{T}(i)}}\Big{)}
\]
\end{enumerate}
\end{enumerate}
\renewcommand{\theenumi}{\arabic{enumi}}
\renewcommand{\labelenumi}{(\theenumi)}
Indeed, (e1) is obvious, and the equations 
$d_{T',i.i-1}=-d_{T,i.i-1},
\,u_{\rho_{T'}(i-1)}=u_{\rho_{T}(i-1)},u_{\rho_{T'}(i)}=u_{\rho_{T}(i)}$
imply (e2)-(c). 
We notice that $u_{\rho_{T}(i-1)}/u_{\rho_{T}(i)}$ equals $1$ if $i$ and $i-1$ 
belong to same Young diagram and $-1$ if not. \par
We shall give the intertwining operator between 
$\pi_{\lambda,\flat}$ and $\pi_{\lambda',\flat}$. 
Using $\psi(T)$ which has given in the equation (1), we define the $K_1$-homomorphism 
$\varPsi_{\flat}$ to be 
\begin{equation*}
\begin{split}
\varPsi_{\flat} : \quad &V_{\lambda} \longrightarrow V_{\lambda'} \\
	&v_T \longmapsto \psi(T)v_{T'} \\
\end{split}
\end{equation*}
We denote the action of $b{\in}\mathscr{H}_{K_1,B_n}(1,q)$ by 
\[ \pi_\lambda(b)v_T=\sum_{S\in \operatorname{STab}(\lambda)}
\Big{(}\pi_\lambda(b)\Big{)}_{T,S}v_{S} \]
where $(\pi_\lambda(b))_{T,S}$'s are elements in $K_1$. 
\begin{Prop}
\begin{equation*}
\begin{split}
\pi_{\lambda'}(b_1)\varPsi_{\flat}(v_{T}) 
	&= \sum_{S\in \operatorname{STab}(\lambda)}
	\Big{(}\pi_\lambda(b_1)\Big{)}_{T,S}\varPsi_{\flat}(v_{S})\\
\pi_{\lambda'}(b_i)\varPsi_{\flat}(v_{T}) 
	&= -\sum_{S\in \operatorname{STab}(\lambda)}
	\Big{(}\pi_\lambda(b_i)\Big{)}_{T,S}\varPsi_{\flat}(v_{S}) \qquad 
	\text{for $i=2,\hdots,n$}  
\end{split}
\end{equation*}
\end{Prop}
\begin{proof}
From (b1) and (e1), it is trivial for $b_1$. 
Assume $i>1$. Comparing (b2) and (e2), the assertion holds for the case 
(a),(b). For the case (c), observing that $u_{\rho_{T}(i-1)}/u_{\rho_{T}(i)}$ equals $1$ 
or $-1$, 
one can prove in the same manner as in the proof of 
Proposition 6.3, just replacing $'^{*}$ with $'$, so we omit the detail. 
\end{proof}
We denote by $(\pi_{\lambda,\flat},V_{\lambda,\flat})$ the restriction of 
$(\pi_\lambda,V_\lambda)$ to $\mathscr{H}_{K_1,B_n}(1,q)^{\flat}$. 
In the same fashions as in Corollary 6.4 and Corollary 6.5, we also obtain the following 
two corollaries. 
\begin{Cor}
$\pi_{\lambda,\flat}{\cong}\pi_{\lambda',\flat}$. 
The $K_1$-homomorphism 
$\varPsi_{\flat}:V_{\lambda,\flat}{\longrightarrow}V_{\lambda',\flat}$ 
is the intertwining operator between $\pi_{\lambda,\flat}$ and 
$\pi_{\lambda',\flat}$. 
\end{Cor}
\begin{Cor}
Let $b{\in}\mathscr{H}_{K_1,B_n}(1,q)^{\flat}$. Then
\[ \Big{(}\pi_{\lambda'}(b)\Big{)}_{T',S'}\psi(T)
	= \Big{(}\pi_{\lambda}(b)\Big{)}_{T,S}\psi(S)
\]
\end{Cor}
Henceforth, we consider $(\pi_{\lambda},V_{\lambda})$ over $\Bar{K}_1$ until the end of 
this section. Assume that $\lambda'=\lambda$. 
We define $\operatorname{STab(\lambda)}^{\flat}$ and $\operatorname{STab(\lambda)}^{-\flat}$ 
to be the sets of all the standard tableaux 
$T$ of shape $\lambda$ such that the smallest number which is not assigned in diagonals of 
Young diagrams appears in the first row and the first column respectively. 
Both $\operatorname{STab(\lambda)}^{\flat}$ and $\operatorname{STab(\lambda)}^{-\flat}$ are 
not empty if $\lambda{\neq}(1,1)$. For $\lambda{\neq}(1,1)$, 
$\operatorname{STab(\lambda)}$ is the disjoint union of 
$\operatorname{STab(\lambda)}^{\flat}$ and $\operatorname{STab(\lambda)}^{-\flat}$. 
There exists an involutive bijection from $\operatorname{STab(\lambda)}^{\flat}$ to 
$\operatorname{STab(\lambda)}^{-\flat}$ which is defined by $T{\mapsto}T'$. 
Now, we define two subspaces $V_{\lambda,\flat}^+$ and 
$V_{\lambda,\flat}^-$ of $V_\lambda$ as follows. 
\begin{equation*}
\begin{split}
V_{\lambda,\flat}^+ &= \bigoplus_{T{\in}\operatorname{STab}(\lambda)^{\flat}}
\Bar{K}_1(\sqrt{\psi(T')}v_{T}+\sqrt{\psi(T)}v_{T'}) \\
V_{\lambda,\flat}^- &= \bigoplus_{T{\in}\operatorname{STab}(\lambda)^{\flat}}
\Bar{K}_1(\sqrt{\psi(T')}v_{T}-\sqrt{\psi(T)}v_{T'})
\end{split}
\end{equation*}
Then, it is clear that $\dim_{\Bar{K}_1}V_{\lambda,\flat}^+ = 
\dim_{\Bar{K}_1}V_{\lambda,\flat}^- = \dim_{\Bar{K}_1}V_{\lambda,\flat}/2$ and 
$V_{\lambda,\flat}$ is a direct sum of $V_{\lambda,\flat}^+$ and 
$V_{\lambda,\flat}^-$ as vector spaces over $\Bar{K}_1$.
\begin{Prop}
If $\lambda'=\lambda$ and $\lambda{\neq}(1,1)$, then 
$V_{\lambda,\flat}=V_{\lambda,\flat}^+{\oplus}V_{\lambda,\flat}^-$ is a 
direct sum decomposition as $\mathscr{H}_{\Bar{K}_1,B_n}(1,q)^{\flat}$-submodules. 
\end{Prop}
\begin{proof}
For $T{\in}\operatorname{STab}(\lambda)^{\flat}$, we have from Corollary 6.11 
\begin{equation*}
\begin{split}
&\pi_{\lambda,\flat}(b)\Big{(}\sqrt{\psi(T')}v_{T}{\pm}\sqrt{\psi(T)}v_{T'}\Big{)} \\
	&= \sum_{S{\in}\operatorname{STab}(\lambda)}\Big{\{}
		\Big{(}\pi_{\lambda,\flat}(b)\Big{)}_{T,S}\sqrt{\psi(T')}v_{S}
	 {\pm} \Big{(}\pi_{\lambda,\flat}(b)\Big{)}_{T',S'}
		\sqrt{\psi(T)}v_{S'}\Big{\}} \\
	&= \sum_{S{\in}\operatorname{STab}(\lambda)}
		\sqrt{\Big{(}\pi_{\lambda,\flat}(b)\Big{)}_{T,S}
		\Big{(}\pi_{\lambda,\flat}(b)\Big{)}_{T',S'}}
		\Big{(}\sqrt{\psi(S')}v_{S}
		{\pm}\sqrt{\psi(S)}v_{S'}\Big{)} \\
	&= \sum_{S{\in}\operatorname{STab}(\lambda)^{\flat}}\Big{\{}
		\sqrt{\Big{(}\pi_{\lambda,\flat}(b)\Big{)}_{T,S}
		\Big{(}\pi_{\lambda,\flat}(b)\Big{)}_{T',S'}}
		{\pm}\sqrt{\Big{(}\pi_{\lambda,\flat}(b)\Big{)}_{T,S'} 
		\Big{(}\pi_{\lambda,\flat}(b)\Big{)}_{T',S}}\Big{\}} \\
		&\times\Big{(}\sqrt{\psi(S')}v_{S}
		{\pm}\sqrt{\psi(S)}v_{S'}\Big{)} \\
\end{split}
\end{equation*}
Thus, $\pi_{\lambda,\flat}(b)\Big{(}\sqrt{\psi(T')}v_{T}{\pm}\sqrt{\psi(T)}v_{T'}\Big{)}$ 
is in $V_{\lambda,\flat}^{\pm}$. 
\end{proof}
If $\lambda=(1,1)$, then $T'=T$. Moreover $s_1T=T^*$ for $T{\in}\operatorname{STab}(\lambda)$. 
Since axial distances $d_{T,2,1}=0$ for all $T{\in}\operatorname{STab}(\lambda)$ in this case, 
(b2)-(c) and (e2)-(c) reduce to 
\[
\pi_{\lambda}(b_2){\langle}v_{T},v_{s_1{T}}{\rangle}
={\langle}v_{T},v_{s_1{T}}{\rangle}
\begin{bmatrix}
0 & 1 \\
1 & 0
\end{bmatrix}, {\qquad}
\pi_{\lambda'}(b_2){\langle}v_{T'},v_{s_1{T'}}{\rangle}
={\langle}v_{T'},v_{s_1{T'}}{\rangle} 
\begin{bmatrix}
0 & 1 \\
1 & 0
\end{bmatrix}.
\]
So letting $T_1=\framebox{1}\; \framebox{2}$, 
we have two submodules 
$V_{\lambda,\flat}^+=\Bar{K}_1v_{T_1},V_{\lambda,\flat}^-=\Bar{K}_1v_{T^{*}_1}$. 
One can readily see that 
$V_{\lambda,\flat}=V_{\lambda,\flat}^+{\oplus}V_{\lambda,\flat}^-$ as 
$\mathscr{H}_{\Bar{K}_1,B_n}(1,q)^{\flat}$-submodules. 
\begin{Theo}
\begin{enumerate}
\rm \item \it
$\pi_{\lambda',\flat}{\cong}\pi_{\lambda,\flat}$.
\rm \item \it
If $\lambda'={\lambda}$, then $\pi_{\lambda,\flat}$ decomposes into 
inequivalent subrepresentations 
$\pi_{\lambda,\flat}^+$ and $\pi_{\lambda,\flat}^-$ over $\Bar{K}_1$, 
the algebraic closure of $K_1$. Furthermore, $\deg\pi_{\lambda,\flat}^+
=\deg\pi_{\lambda,\flat}^-=\deg\pi_{\lambda,\flat}/2$.
\rm \item \it
Let $\{\lambda_i,\lambda_i',\mu_j\}_{i,j}$ 
be the set of $2$-tuples of Young diagrams 
such that 
\begin{equation*}
\lambda_i'{\neq}\lambda_i,{\quad}\mu_j'=\mu_j
\end{equation*}
Then 
\[
\operatorname{Irr}(\mathscr{H}_{\Bar{K}_1,B_n}(q)^{\flat})=
\{\pi_{\lambda_i,\flat},\pi_{\mu_j,\flat}^+,\pi_{\mu_j,\flat}^-,\}_{i,j}
\]
is a basic set of irreducible representations of $\mathscr{H}_{\Bar{K}_1,B_n}(q)^{\flat}$. 
$\mathscr{H}_{\Bar{K}_1,B_n}(q)^{\flat}$ is semisimple. 
\end{enumerate}
\end{Theo}
\begin{proof}
(1) is Corollary 6.10 itself. 
The submodules $V_{\lambda,\flat}^+$ and $V_{\lambda,\flat}^-$ 
correspond to the case (2) of Theorem 2.7. Hence they are simple and mutually 
non-isomorphic. Thus (2) holds. 
The proof of (3) is in the same manner as in the proof of Theorem 6.2, 
just replacing $\natural$ and $\bar{K}_0$ with $\flat$ and $\bar{K}_1$ respectively. 
So we omit the detail. 
\end{proof}
\section{The branching rule for the Hecke algebras of Type $D$}
Henceforth we assume $n{\geq}4$. 
We denote by $(\pi_{\lambda,\dagger},V_{\lambda,\dagger})$, 
$(\pi_{\lambda,\natural,\dagger},V_{\lambda,\natural,\dagger})$, 
$(\pi_{\lambda,\sharp,\dagger},V_{\lambda,\sharp,\dagger})$, 
$(\pi_{\lambda,\flat,\dagger},V_{\lambda,\flat,\dagger})$ the restrictions of 
$(\pi_\lambda,V_\lambda)$,$(\pi_{\lambda,\natural},V_{\lambda,\natural})$,
$(\pi_{\lambda,\sharp},V_{\lambda,\sharp})$,
$(\pi_{\lambda,\flat},V_{\lambda,\flat})$ to $\mathscr{H}_{K_1,B_n}(1,q)^{\dagger}$ respectively. 
It is clear that 
$\pi_{\lambda,\dagger}{\cong}\pi_{\lambda,\natural,\dagger}{\cong}
\pi_{\lambda,\sharp,\dagger}{\cong}\pi_{\lambda,\flat,\dagger}$. 
If $\lambda^*=\lambda$, then we denote by 
$(\pi_{\lambda,\sharp,\dagger}^{\pm},V_{\lambda,\sharp,\dagger}^{\pm})$ the restriction of 
$(\pi_{\lambda,\sharp}^{\pm},V_{\lambda,\sharp}^{\pm})$ to $\mathscr{H}_{K_1,B_n}(1,q)^{\dagger}$. 
We adopt the similar notations for the cases $\lambda'=\lambda$ and $\lambda'^*=\lambda$. 
We notice that $\mathscr{H}_{K_1,B_n}(1,q)^{\sharp}=\mathscr{H}_{K_1,D_n}(q)$ and 
$\mathscr{H}_{K_1,B_n}(1,q)^{\flat}$ and $\mathscr{H}_{K_1,B_n}(1,q)^{\natural}$ are 
all $\mathbb{Z}_2$-crossed products with $A_{\bar{0}}=\mathscr{H}_{K_1,B_n}(1,q)^{\dagger}$
\begin{Prop}
$\pi_{\lambda,\dagger}{\cong}\pi_{\lambda^{*},\dagger}{\cong}
\pi_{\lambda',\dagger}{\cong}\pi_{\lambda'^*,\dagger}$
\end{Prop}
\begin{proof}
Since $\mathscr{H}_{K_1,B_n}(1,q)^{\dagger}=\mathscr{H}_{K_1,B_n}(1,q)^{\sharp}{\cap}
\mathscr{H}_{K_1,B_n}(1,q)^{\flat}{\cap}\mathscr{H}_{K_1,B_n}(1,q)^{\natural}$, 
we have that $\pi_{\lambda,\dagger},\pi_{\lambda^*,\dagger},\pi_{\lambda',\dagger}$ and 
$\pi_{\lambda'^*,\dagger}$ are all equivalent. 
\end{proof}
All the representations $V_{\lambda}$,$V_{\lambda,\sharp}$,$V_{\lambda,\flat}$,
$V_{\lambda,\natural}$,${\cdots}$ can be defined over 
$\Bar{K}_1$, so we consider all these over $\bar{K}_1$. 
At first, we assume that $\lambda$, $\lambda^*$, $\lambda'$ and $\lambda'^*$ are mutually 
different. 
\begin{Prop}
Let $\lambda$, $\lambda^*$, $\lambda'$ and $\lambda'^*$ are mutually different. 
Then $\pi_{\lambda,\sharp,\dagger}$ is irreducible. 
\end{Prop}
\begin{proof}
We have already shown that $\pi_{\lambda,\sharp}$ and $\pi_{\lambda',\sharp}$ are 
irreducible and inequivalent. 
On the other hand, $\pi_{\lambda,\sharp,\dagger}{\cong}\pi_{\lambda',\sharp,\dagger}$ 
by Proposition 7.1. This corresponds to the case (1) of Theorem 2.7. 
Thus we can conclude that $\pi_{\lambda,\sharp,\dagger}$ is irreducible. 
\end{proof}
Next, we assume that only one of the three cases $\lambda^*=\lambda$, 
$\lambda'=\lambda$, $\lambda'^*=\lambda$ holds, and considering representations 
over $\Bar{K}_1$, the algebraic closure of $K_1$
\begin{Prop}
\noindent
\begin{enumerate}
\rm\item\it
If $\lambda^*=\lambda$, $\lambda'{\;\neq\;}\lambda$, $\lambda'^*{\;\neq\;}\lambda$, 
then $\pi_{\lambda,\sharp,\dagger}^+{\cong}\pi_{\lambda',\sharp,\dagger}^+$ and 
$\pi_{\lambda,\sharp,\dagger}^-{\cong}\pi_{\lambda',\sharp,\dagger}^-$. 
They are irreducible and mutually inequivalent. 
\rm\item\it
If $\lambda^*{\;\neq\;}\lambda$, $\lambda'=\lambda$, $\lambda'^*{\;\neq\;}\lambda$, 
then $\pi_{\lambda,\flat,\dagger}^+{\cong}\pi_{\lambda',\flat,\dagger}^+$ and 
$\pi_{\lambda,\flat,\dagger}^-{\cong}\pi_{\lambda',\flat,\dagger}^-$. 
They are irreducible and mutually inequivalent. 
\rm\item\it
If $\lambda^*{\;\neq\;}\lambda$, $\lambda'{\;\neq\;}\lambda$, $\lambda'^*=\lambda$, 
then $\pi_{\lambda,\natural,\dagger}^+{\cong}\pi_{\lambda',\natural,\dagger}^+$ and 
$\pi_{\lambda,\natural,\dagger}^-{\cong}\pi_{\lambda',\natural,\dagger}^-$. 
They are irreducible and mutually inequivalent. 
\end{enumerate}
\end{Prop}
\begin{proof}
When $\lambda^*=\lambda$, $\lambda'^*=\lambda'$ also holds. 
We have already shown that $\pi_{\lambda,\sharp}$ and $\pi_{\lambda',\sharp}$ have the irreducible 
decompositions $\pi_{\lambda,\sharp}=\pi_{\lambda,\sharp}^+{\oplus}\pi_{\lambda,\sharp}^-$ and 
$\pi_{\lambda',\sharp}=\pi_{\lambda',\sharp}^+{\oplus}\pi_{\lambda',\sharp}^-$ respectively. 
$\pi_{\lambda,\sharp}{\ncong}\pi_{\lambda',\sharp}$ since $\lambda'{\;\neq\;}\lambda$, while 
$\pi_{\lambda,\sharp,\dagger}{\cong}\pi_{\lambda',\sharp,\dagger}$. 
Therefore, $V_{\lambda,\sharp,\dagger}^+$ is 
$\mathscr{H}_{\bar{K}_1,B_n}(1,q)^{\dagger}$-isomorphic to a submodule of 
$V_{\lambda',\sharp,\dagger}$. 
This corresponds to the case (1) of Theorem 2.7, hence $\pi_{\lambda,\sharp,\dagger}^+$ is 
irreducible. 
In the same manner, $\pi_{\lambda,\sharp,\dagger}^-$, $\pi_{\lambda',\sharp,\dagger}^+$, 
$\pi_{\lambda',\sharp,\dagger}^-$ are also irreducible. 
Assume $\pi_{\lambda,\sharp,\dagger}^+{\cong}\pi_{\lambda,\sharp,\dagger}^-$. Then four 
irreducible representations $\pi_{\lambda,\sharp,\dagger}^{\pm}$, 
$\pi_{\lambda',\sharp,\dagger}^{\pm}$ are all equivalent. But this is impossible because of 
Theorem 2.7. 
Thus $\pi_{\lambda,\sharp,\dagger}^+{\ncong}\pi_{\lambda,\sharp,\dagger}^-$ holds. 
Consequently, we also have 
$\pi_{\lambda,\sharp,\dagger}^+{\cong}\pi_{\lambda',\sharp,\dagger}^+$ and 
$\pi_{\lambda,\sharp,\dagger}^-{\cong}\pi_{\lambda',\sharp,\dagger}^-$. 
In the same manner, we also obtain (2) and (3). 
\end{proof}
Last, we assume $\lambda^*=\lambda$ and $\lambda'=\lambda$ 
(these imply $\lambda'^*=\lambda$).
For brevity, we set   
\begin{equation*}
\begin{split}
&v_T^{++}=\sqrt[4]{\psi(T')\psi(T'^*)}(v_T+v_{T^*})
+\sqrt[4]{\psi(T)\psi(T^*)}(v_{T'}+v_{T'^*}),\\
&v_T^{+-}=\sqrt[4]{\psi(T')\psi(T'^*)}(v_T+v_{T^*})
-\sqrt[4]{\psi(T)\psi(T^*)}(v_{T'}+v_{T'^*}),\\
&v_T^{-+}=\sqrt[4]{\psi(T')\psi(T'^*)}(v_T-v_{T^*})
+\sqrt[4]{\psi(T)\psi(T^*)}(v_{T'}-v_{T'^*}),\\
&v_T^{--}=\sqrt[4]{\psi(T')\psi(T'^*)}(v_T-v_{T^*})
-\sqrt[4]{\psi(T)\psi(T^*)}(v_{T'}-v_{T'^*}).
\end{split}
\end{equation*}
We define $\Bar{K}_1$-subspaces 
$V_{\lambda,\sharp,\dagger}^{++}$, $V_{\lambda,\sharp,\dagger}^{+-}$, 
$V_{\lambda,\sharp,\dagger}^{-+}$, $V_{\lambda,\sharp,\dagger}^{--}$ of 
$V_{\lambda}$ as follows. 
\begin{equation*}
\begin{split}
V_{\lambda,\sharp,\dagger}^{++}&=\bigoplus_{T\,{\in}\,\operatorname{STab}(\lambda)^{\dagger}}
\Bar{K}_1v_T^{++},{\qquad}
V_{\lambda,\sharp,\dagger}^{+-}=\bigoplus_{T\,{\in}\,\operatorname{STab}(\lambda)^{\dagger}}
\Bar{K}_1v_T^{+-},\\
V_{\lambda,\sharp,\dagger}^{-+}&=\bigoplus_{T\,{\in}\,\operatorname{STab}(\lambda)^{\dagger}}
\Bar{K}_1v_T^{-+},{\qquad}
V_{\lambda,\sharp,\dagger}^{--}=\bigoplus_{T\,{\in}\,\operatorname{STab}(\lambda)^{\dagger}}
\Bar{K}_1v_T^{--}
\end{split}
\end{equation*}
where $\operatorname{STab}(\lambda)^{\dagger}=\operatorname{STab}(\lambda)^{\sharp}
{\cap}\operatorname{STab}(\lambda)^{\flat}$. 
\begin{Rem}
We notice that each fourth square root has four branches. 
We take a suitable branch in each computation of fourth square roots for consistency. 
\end{Rem}
\begin{Prop}
\[
V_{\lambda,\sharp,\dagger}^{++}{\subseteq}V_{\lambda,\sharp}^+{\cap}V_{\lambda,\flat}^+, {\quad}
V_{\lambda,\sharp,\dagger}^{+-}{\subseteq}V_{\lambda,\sharp}^+{\cap}V_{\lambda,\flat}^-, {\quad}
V_{\lambda,\sharp,\dagger}^{-+}{\subseteq}V_{\lambda,\sharp}^-{\cap}V_{\lambda,\flat}^+, {\quad}
V_{\lambda,\sharp,\dagger}^{--}{\subseteq}V_{\lambda,\sharp}^-{\cap}V_{\lambda,\flat}^-,
\]
and furthermore, we have direct sum decompositions as 
$\Bar{K}_1$-vector spaces as follows. 
\begin{equation*}
\begin{split}
V_{\lambda,\sharp}^+&=V_{\lambda,\sharp,\dagger}^{++}{\oplus}V_{\lambda,\sharp,\dagger}^{+-},{\quad}
V_{\lambda,\sharp}^-=V_{\lambda,\sharp,\dagger}^{-+}{\oplus}V_{\lambda,\sharp,\dagger}^{--},\\
V_{\lambda,\flat}^+&=V_{\lambda,\sharp,\dagger}^{++}{\oplus}V_{\lambda,\sharp,\dagger}^{-+},{\quad}
V_{\lambda,\flat}^-=V_{\lambda,\sharp,\dagger}^{+-}{\oplus}V_{\lambda,\sharp,\dagger}^{--},\\
V_{\lambda,\natural}^+&=V_{\lambda,\sharp,\dagger}^{++}{\oplus}V_{\lambda,\sharp,\dagger}^{--},{\quad}
V_{\lambda,\natural}^-=V_{\lambda,\sharp,\dagger}^{+-}{\oplus}V_{\lambda,\sharp,\dagger}^{-+}.
\end{split}
\end{equation*}
\end{Prop}
\begin{proof}
Assume $T{\in}\operatorname{STab}(\lambda)^{\dagger}$. Then $T$ and $T'$ belong to 
$\operatorname{STab}(\lambda)^{\sharp}$. Therefore, we have 
$V_{\lambda,\sharp,\dagger}^{++},V_{\lambda,\sharp,\dagger}^{+-}
{\in}V_{\lambda,\sharp}^+$ and 
$V_{\lambda,\sharp,\dagger}^{-+},V_{\lambda,\sharp,\dagger}^{--}
{\in}V_{\lambda,\sharp}^-$ immediately. 
We also see that $T,T^*{\in}\operatorname{STab}(\lambda)^{\flat}$ and 
$T,T'{\in}\operatorname{STab}(\lambda)^{\natural}$. 
If $n/2$ is even, then $\psi(T^*)=\psi(T)$ and hence we have two expressions for each linear 
combination as follows. 
\begin{equation*}
\begin{split}
&v_T^{++}=
\begin{cases}
\sqrt{\psi(T')}v_T+\sqrt{\psi(T)}v_{T'}
+\sqrt{\psi(T'^*)}v_{T^*}+\sqrt{\psi(T^*)}v_{T'^*}
{\in}V_{\lambda,\flat}^+,\\
\sqrt{\psi(T'^*)}v_T+\sqrt{\psi(T)}v_{T'^*}
+\sqrt{\psi(T^*)}v_{T'}+\sqrt{\psi(T')}v_{T^*}
{\in}V_{\lambda,\natural}^+,
\end{cases}\\
&v_T^{+-}=
\begin{cases}
\sqrt{\psi(T')}v_T-\sqrt{\psi(T)}v_{T'}
+\sqrt{\psi(T'^*)}v_{T^*}-\sqrt{\psi(T^*)}v_{T'^*}
{\in}V_{\lambda,\flat}^-,\\
\sqrt{\psi(T'^*)}v_T-\sqrt{\psi(T)}v_{T'^*}
-(\sqrt{\psi(T^*)}v_{T'}-\sqrt{\psi(T')}v_{T^*})
{\in}V_{\lambda,\natural}^-,
\end{cases}\\
&v_T^{-+}=
\begin{cases}
\sqrt{\psi(T')}v_T+\sqrt{\psi(T)}v_{T'}
-(\sqrt{\psi(T'^*)}v_{T^*}+\sqrt{\psi(T^*)}v_{T'^*})
{\in}V_{\lambda,\flat}^+,\\
\sqrt{\psi(T'^*)}v_T-\sqrt{\psi(T)}v_{T'^*}
+\sqrt{\psi(T^*)}v_{T'}-\sqrt{\psi(T')}v_{T^*}
{\in}V_{\lambda,\natural}^-,
\end{cases}\\
&v_T^{--}=
\begin{cases}
\sqrt{\psi(T')}v_T-\sqrt{\psi(T)}v_{T'}
-(\sqrt{\psi(T'^*)}v_{T^*}-\sqrt{\psi(T^*)}v_{T'^*})
{\in}V_{\lambda,\flat}^-,\\
\sqrt{\psi(T'^*)}v_T+\sqrt{\psi(T)}v_{T'^*}
-(\sqrt{\psi(T^*)}v_{T'}+\sqrt{\psi(T')}v_{T^*})
{\in}V_{\lambda,\natural}^+.
\end{cases}
\end{split}
\end{equation*}
If $n/2$ is odd, then $\psi(T^*)=-\psi(T)$ and hence 
\begin{equation*}
\begin{split}
&v_T^{++}=
\begin{cases}
\sqrt[4]{-1}\big{\{}\sqrt{\psi(T')}v_T+\sqrt{\psi(T)}v_{T'}
+\sqrt{\psi(T'^*)}v_{T^*}+\sqrt{\psi(T^*)}v_{T'^*}\big{\}}
{\in}V_{\lambda,\flat}^+,\\
\sqrt[4]{-1}\big{\{}\sqrt{\psi(T'^*)}v_T+\sqrt{\psi(T)}v_{T'^*}
+\sqrt{\psi(T^*)}v_{T'}+\sqrt{\psi(T')}v_{T^*}\big{\}}
{\in}V_{\lambda,\natural}^+,
\end{cases}\\
&v_T^{+-}=
\begin{cases}
\sqrt[4]{-1}\big{\{}\sqrt{\psi(T')}v_T-\sqrt{\psi(T)}v_{T'}
+\sqrt{\psi(T'^*)}v_{T^*}-\sqrt{\psi(T^*)}v_{T'^*}\big{\}}
{\in}V_{\lambda,\flat}^-,\\
\sqrt[4]{-1}\big{\{}\sqrt{\psi(T'^*)}v_T-\sqrt{\psi(T)}v_{T'^*}
-(\sqrt{\psi(T^*)}v_{T'}-\sqrt{\psi(T')}v_{T^*})\big{\}}
{\in}V_{\lambda,\natural}^-,
\end{cases}\\
&v_T^{-+}=
\begin{cases}
\sqrt[4]{-1}\big{\{}\sqrt{\psi(T')}v_T+\sqrt{\psi(T)}v_{T'}
-(\sqrt{\psi(T'^*)}v_{T^*}+\sqrt{\psi(T^*)}v_{T'^*})\big{\}}
{\in}V_{\lambda,\flat}^+,\\
\sqrt[4]{-1}\big{\{}\sqrt{\psi(T'^*)}v_T-\sqrt{\psi(T)}v_{T'^*}
+\sqrt{\psi(T^*)}v_{T'}-\sqrt{\psi(T')}v_{T^*}\big{\}}
{\in}V_{\lambda,\natural}^-,
\end{cases}\\
&v_T^{--}=
\begin{cases}
\sqrt[4]{-1}\big{\{}\sqrt{\psi(T')}v_T-\sqrt{\psi(T)}v_{T'}
-(\sqrt{\psi(T'^*)}v_{T^*}-\sqrt{\psi(T^*)}v_{T'^*})\big{\}}
{\in}V_{\lambda,\flat}^-,\\
\sqrt[4]{-1}\big{\{}\sqrt{\psi(T'^*)}v_T+\sqrt{\psi(T)}v_{T'^*}
-(\sqrt{\psi(T^*)}v_{T'}+\sqrt{\psi(T')}v_{T^*})\big{\}}
{\in}V_{\lambda,\natural}^+.\\
\end{cases}
\end{split}
\end{equation*}
Thus the first assertion holds. 
We also obtain that 
\begin{equation*}
\begin{split}
&v_T^{++}+v_T^{+-}=2\sqrt[4]{\psi(T')\psi(T'^*)}(v_T+v_{T^*}),\\
&v_T^{++}-v_T^{+-}=2\sqrt[4]{\psi(T)\psi(T^*)}(v_{T'}+v_{T'^*}).\\
\end{split}
\end{equation*}
Thus every $v_T+v_{T^*}$ ($T{\in}\operatorname{STab}(\lambda)^{\sharp}$) belongs to 
$V_{\lambda,\sharp,\dagger}^{++}+V_{\lambda,\sharp,\dagger}^{+-}$. 
We readily see that $V_{\lambda,\sharp,\dagger}^{++}+V_{\lambda,\sharp,\dagger}^{+-}$ is 
a direct sum by the first assertion. 
Therefore we obtain 
$V_{\lambda,\sharp}^+=V_{\lambda,\sharp,\dagger}^{++}{\oplus}V_{\lambda,\sharp,\dagger}^{+-}$.
The others are obtained in the same manner. 
\end{proof}
\begin{Prop}
$V_{\lambda,\sharp,\dagger}^{++},V_{\lambda,\sharp,\dagger}^{+-},
V_{\lambda,\sharp,\dagger}^{-+},V_{\lambda,\sharp,\dagger}^{--}$ are 
$\mathscr{H}_{\Bar{K}_1,B_n}(1,q)^{\dagger}$-modules. 
\end{Prop}
\begin{proof}
Using Corollary 6.5 and Corollary 6.11, one can check by a direct computation that  
\begin{equation*}
\begin{split}
&\pi_{\lambda}(b)\big{(}\sqrt[4]{\psi(T')\psi(T'^*)}v_T\big{)}\\
&=\sum_{S{\in}\operatorname{STab}(\lambda)^{\dagger}}
\Big{\{}\sqrt[4]{\big{(}\pi_{\lambda}(b)\big{)}_{TS}\big{(}\pi_{\lambda}(b)\big{)}_{T^*S^*}
\big{(}\pi_{\lambda}(b)\big{)}_{T'S'}\big{(}\pi_{\lambda}(b)\big{)}_{T'^*S'^*}}
\sqrt[4]{\psi(S')\psi(S'^*)}v_{S}\\
&+\sqrt[4]{\big{(}\pi_{\lambda}(b)\big{)}_{TS'}\big{(}\pi_{\lambda}(b)\big{)}_{T^*S'^*}
\big{(}\pi_{\lambda}(b)\big{)}_{T'S}\big{(}\pi_{\lambda}(b)\big{)}_{T'^*S^*}}
\sqrt[4]{\psi(S)\psi(S^*)}v_{S'}\\
&+\sqrt[4]{\big{(}\pi_{\lambda}(b)\big{)}_{TS^*}\big{(}\pi_{\lambda}(b)\big{)}_{T^*S}
\big{(}\pi_{\lambda}(b)\big{)}_{T'S'^*}\big{(}\pi_{\lambda}(b)\big{)}_{T'^*S'}}
\sqrt[4]{\psi(S'^*)\psi(S')}v_{S^*}\\
&+\sqrt[4]{\big{(}\pi_{\lambda}(b)\big{)}_{TS'^*}\big{(}\pi_{\lambda}(b)\big{)}_{T^*S{'}}
\big{(}\pi_{\lambda}(b)\big{)}_{T'S^*}\big{(}\pi_{\lambda}(b)\big{)}_{T'^*S}}
\sqrt[4]{\psi(S^*)\psi(S)}v_{S'^*}\Big{\}}
\end{split}
\end{equation*}
for $b{\in}\mathscr{H}_{\Bar{K}_1,B_n}(1,q)^{\dagger}$. 
Similarly we have 
\begin{equation*}
\begin{split}
&\pi_{\lambda}(b)\big{(}\sqrt[4]{\psi(T)\psi(T^*)}v_{T'}\big{)}\\
&=\sum_{S{\in}\operatorname{STab}(\lambda)^{\dagger}}
\Big{\{}\sqrt[4]{\big{(}\pi_{\lambda}(b)\big{)}_{T'S}\big{(}\pi_{\lambda}(b)\big{)}_{T'^*S^*}
\big{(}\pi_{\lambda}(b)\big{)}_{TS'}\big{(}\pi_{\lambda}(b)\big{)}_{T^*S'^*}}
\sqrt[4]{\psi(S')\psi(S'^*)}v_{S}\\
&+\sqrt[4]{\big{(}\pi_{\lambda}(b)\big{)}_{T'S'}\big{(}\pi_{\lambda}(b)\big{)}_{T'^*S'^*}
\big{(}\pi_{\lambda}(b)\big{)}_{TS}\big{(}\pi_{\lambda}(b)\big{)}_{T^*S^*}}
\sqrt[4]{\psi(S)\psi(S^*)}v_{S'}\\
&+\sqrt[4]{\big{(}\pi_{\lambda}(b)\big{)}_{T'S^*}\big{(}\pi_{\lambda}(b)\big{)}_{T'^*S}
\big{(}\pi_{\lambda}(b)\big{)}_{TS'^*}\big{(}\pi_{\lambda}(b)\big{)}_{T^*S'}}
\sqrt[4]{\psi(S'^*)\psi(S')}v_{S^*}\\
&+\sqrt[4]{\big{(}\pi_{\lambda}(b)\big{)}_{T'S'^*}\big{(}\pi_{\lambda}(b)\big{)}_{T'^*S'}
\big{(}\pi_{\lambda}(b)\big{)}_{TS^*}\big{(}\pi_{\lambda}(b)\big{)}_{T^*S}}
\sqrt[4]{\psi(S^*)\psi(S)}v_{S'^*}\Big{\}}
\end{split}
\end{equation*}
\begin{equation*}
\begin{split}
&\pi_{\lambda}(b)\big{(}\sqrt[4]{\psi(T'^*)\psi(T')}v_{T^*}\big{)}\\
&=\sum_{S{\in}\operatorname{STab}(\lambda)^{\dagger}}
\Big{\{}\sqrt[4]{\big{(}\pi_{\lambda}(b)\big{)}_{T^*S}\big{(}\pi_{\lambda}(b)\big{)}_{TS^*}
\big{(}\pi_{\lambda}(b)\big{)}_{T'^*S'}\big{(}\pi_{\lambda}(b)\big{)}_{T'S'^*}}
\sqrt[4]{\psi(S')\psi(S'^*)}v_{S}\\
&+\sqrt[4]{\big{(}\pi_{\lambda}(b)\big{)}_{T^*S'}\big{(}\pi_{\lambda}(b)\big{)}_{TS'^*}
\big{(}\pi_{\lambda}(b)\big{)}_{T'^*S}\big{(}\pi_{\lambda}(b)\big{)}_{T'S^*}}
\sqrt[4]{\psi(S)\psi(S^*)}v_{S'}\\
&+\sqrt[4]{\big{(}\pi_{\lambda}(b)\big{)}_{T^*S^*}\big{(}\pi_{\lambda}(b)\big{)}_{TS}
\big{(}\pi_{\lambda}(b)\big{)}_{T'^*S'^*}\big{(}\pi_{\lambda}(b)\big{)}_{T'S'}}
\sqrt[4]{\psi(S'^*)\psi(S')}v_{S^*}\\
&+\sqrt[4]{\big{(}\pi_{\lambda}(b)\big{)}_{T^*S'^*}\big{(}\pi_{\lambda}(b)\big{)}_{TS'}
\big{(}\pi_{\lambda}(b)\big{)}_{T'^*S^*}\big{(}\pi_{\lambda}(b)\big{)}_{T'S}}
\sqrt[4]{\psi(S^*)\psi(S)}v_{S'^*}\Big{\}}
\end{split}
\end{equation*}
\begin{equation*}
\begin{split}
&\pi_{\lambda}(b)\big{(}\sqrt[4]{\psi(T')\psi(T'^*)}v_{T'^*}\big{)}\\
&=\sum_{S{\in}\operatorname{STab}(\lambda)^{\dagger}}
\Big{\{}\sqrt[4]{\big{(}\pi_{\lambda}(b)\big{)}_{T'^*S}\big{(}\pi_{\lambda}(b)\big{)}_{T'S^*}
\big{(}\pi_{\lambda}(b)\big{)}_{T^*S'}\big{(}\pi_{\lambda}(b)\big{)}_{TS'^*}}
\sqrt[4]{\psi(S')\psi(S'^*)}v_{S}\\
&+\sqrt[4]{\big{(}\pi_{\lambda}(b)\big{)}_{T'^*S'}\big{(}\pi_{\lambda}(b)\big{)}_{T'S'^*}
\big{(}\pi_{\lambda}(b)\big{)}_{T^*S}\big{(}\pi_{\lambda}(b)\big{)}_{TS^*}}
\sqrt[4]{\psi(S)\psi(S^*)}v_{S'}\\
&+\sqrt[4]{\big{(}\pi_{\lambda}(b)\big{)}_{T'^*S^*}\big{(}\pi_{\lambda}(b)\big{)}_{T'S}
\big{(}\pi_{\lambda}(b)\big{)}_{T^*S'^*}\big{(}\pi_{\lambda}(b)\big{)}_{TS'}}
\sqrt[4]{\psi(S'^*)\psi(S')}v_{S^*}\\
&+\sqrt[4]{\big{(}\pi_{\lambda}(b)\big{)}_{T'^*S'^*}\big{(}\pi_{\lambda}(b)\big{)}_{T'S'}
\big{(}\pi_{\lambda}(b)\big{)}_{T^*S^*}\big{(}\pi_{\lambda}(b)\big{)}_{TS}}
\sqrt[4]{\psi(S^*)\psi(S)}v_{S'^*}\Big{\}}
\end{split}
\end{equation*}
Summing up all of those, we get 
\begin{equation*}
\begin{split}
&\pi_{\lambda}(b)v_T^{++}
=\sum_{S{\in}\operatorname{STab}(\lambda)^{\dagger}}
\Big{\{}\sqrt[4]{\big{(}\pi_{\lambda}(b)\big{)}_{TS}\big{(}\pi_{\lambda}(b)\big{)}_{T^*S^*}
\big{(}\pi_{\lambda}(b)\big{)}_{T'S'}\big{(}\pi_{\lambda}(b)\big{)}_{T'^*S'^*}}\\
&+\sqrt[4]{\big{(}\pi_{\lambda}(b)\big{)}_{T'S}\big{(}\pi_{\lambda}(b)\big{)}_{T'^*S^*}
\big{(}\pi_{\lambda}(b)\big{)}_{TS'}\big{(}\pi_{\lambda}(b)\big{)}_{T^*S'^*}}\\
&+\sqrt[4]{\big{(}\pi_{\lambda}(b)\big{)}_{T^*S}\big{(}\pi_{\lambda}(b)\big{)}_{TS^*}
\big{(}\pi_{\lambda}(b)\big{)}_{T'^*S'}\big{(}\pi_{\lambda}(b)\big{)}_{T'S'^*}}\\
&+\sqrt[4]{\big{(}\pi_{\lambda}(b)\big{)}_{T'^*S}\big{(}\pi_{\lambda}(b)\big{)}_{T'S^*}
\big{(}\pi_{\lambda}(b)\big{)}_{T^*S'}\big{(}\pi_{\lambda}(b)\big{)}_{TS'^*}}\Big{\}}
v_S^{++}.
\end{split}
\end{equation*}
Thus $\pi_{\lambda}(b)V_{\lambda,\sharp,\dagger}^{++}{\subseteq}V_{\lambda,\sharp,\dagger}^{++}$ 
for $b{\in}\mathscr{H}_{\Bar{K}_1,B_n}(1,q)^{\dagger}$. 
In the same fashion, we also obtain 
$\pi_{\lambda}(b)V_{\lambda,\sharp,\dagger}^{+-}{\subseteq}V_{\lambda,\sharp,\dagger}^{+-}$,
$\pi_{\lambda}(b)V_{\lambda,\sharp,\dagger}^{-+}{\subseteq}V_{\lambda,\sharp,\dagger}^{-+}$ and 
$\pi_{\lambda}(b)V_{\lambda,\sharp,\dagger}^{--}{\subseteq}V_{\lambda,\sharp,\dagger}^{--}$.
\end{proof}
We denote the representations of $\mathscr{H}_{\Bar{K}_1,B_n}(1,q)^{\dagger}$ given 
in Proposition 7.6 by $(\pi_{\lambda,\sharp,\dagger}^{++},V_{\lambda,\sharp,\dagger}^{++})$, 
$(\pi_{\lambda,\sharp,\dagger}^{+-},V_{\lambda,\sharp,\dagger}^{+-})$, 
$(\pi_{\lambda,\sharp,\dagger}^{-+},V_{\lambda,\sharp,\dagger}^{-+})$ and 
$(\pi_{\lambda,\sharp,\dagger}^{--},V_{\lambda,\sharp,\dagger}^{--})$. 
\begin{Prop}
$(\pi_{\lambda,\sharp,\dagger}^{++},V_{\lambda,\sharp,\dagger}^{++})$, 
$(\pi_{\lambda,\sharp,\dagger}^{+-},V_{\lambda,\sharp,\dagger}^{+-})$, 
$(\pi_{\lambda,\sharp,\dagger}^{-+},V_{\lambda,\sharp,\dagger}^{-+})$, 
$(\pi_{\lambda,\sharp,\dagger}^{--},V_{\lambda,\sharp,\dagger}^{--})$ 
are mutually non-isomorphic simple left $\mathscr{H}_{\Bar{K}_1,B_n}(1,q)^{\dagger}$-modules. 
\end{Prop}
\begin{proof}
Observing that $\mathscr{H}_{\Bar{K}_1,B_n}(1,q)^{\sharp}$ and 
$\mathscr{H}_{\Bar{K}_1,B_n}(1,q)^{\flat}$ and $\mathscr{H}_{\Bar{K}_1,B_n}(1,q)^{\natural}$ are 
all $\mathbb{Z}_2$-crossed products with $A_{\bar{0}}=\mathscr{H}_{\Bar{K}_1,B_n}(1,q)^{\dagger}$, 
one can deduce the assertion directly from Proposition 7.5, Proposition 7.6 and Theorem 2.7 (2). 
\end{proof}
We have the following immediately from Proposition 7.5, Proposition 7.6 and Proposition 7.7. 
\begin{Cor}
We have branching rules for $\lambda$ such that $\lambda=\lambda^*=\lambda'=\lambda'^*$ 
as follows. 
\begin{equation*}
\begin{split}
\pi_{\lambda,\dagger}
&=\pi_{\lambda,\sharp,\dagger}^{++}{\oplus}\pi_{\lambda,\sharp,\dagger}^{+-}{\oplus}
\pi_{\lambda,\sharp,\dagger}^{-+}{\oplus}\pi_{\lambda,\sharp,\dagger}^{--},\\
\pi_{\lambda,\sharp,\dagger}^+
&=\pi_{\lambda,\sharp,\dagger}^{++}{\oplus}\pi_{\lambda,\sharp,\dagger}^{+-},{\qquad}
\pi_{\lambda,\sharp,\dagger}^-
=\pi_{\lambda,\sharp,\dagger}^{-+}{\oplus}\pi_{\lambda,\sharp,\dagger}^{--},\\
\pi_{\lambda,\flat,\dagger}^+
&=\pi_{\lambda,\sharp,\dagger}^{++}{\oplus}\pi_{\lambda,\sharp,\dagger}^{-+},{\qquad}
\pi_{\lambda,\flat,\dagger}^-
=\pi_{\lambda,\sharp,\dagger}^{+-}{\oplus}\pi_{\lambda,\sharp,\dagger}^{--},\\
\pi_{\lambda,\natural,\dagger}^+
&=\pi_{\lambda,\sharp,\dagger}^{++}{\oplus}\pi_{\lambda,\sharp,\dagger}^{--},{\qquad}
\pi_{\lambda,\natural,\dagger}^-
=\pi_{\lambda,\sharp,\dagger}^{+-}{\oplus}\pi_{\lambda,\sharp,\dagger}^{-+},\\
\end{split}
\end{equation*}
$\pi_{\lambda,\sharp,\dagger}^{++},\pi_{\lambda,\sharp,\dagger}^{+-},
\pi_{\lambda,\sharp,\dagger}^{-+},\pi_{\lambda,\sharp,\dagger}^{--}$ are irreducible and 
mutually inequivalent. They have the same degree. 
\end{Cor}
\begin{Theo}
Let $\{\lambda_i,\lambda_i^*,\lambda_i',\lambda_i'^*,
\mu_j,\mu_j{'},\nu_k,\nu_k^*,\kappa_l,\kappa_l',
\iota_m\}_{i,j,k,l,m}$ 
be the set of $2$-tuples of Young diagrams 
such that 
\begin{equation*}
\begin{split}
&\lambda_i{'}{\;\neq\;}\lambda_i,\lambda_i^*{\;\neq\;}\lambda_i,\lambda_i'^*{\;\neq\;}\lambda_i,{\quad}
\mu_j^*=\mu_j,\mu_j{'}{\;\neq\;}\mu_j,\mu_j'^*{\;\neq\;}\mu_j,\\
&\nu_k^*{\;\neq\;}\nu_k,\nu_k{'}=\nu_k,\nu_k'^*{\;\neq\;}\nu_k,{\quad}
\kappa_l^*{\;\neq\;}\kappa_l,\kappa_l{'}{\;\neq\;}\kappa_l,\kappa_l'^*=\kappa_l,{\quad}
\iota_m{'}=\iota_m^*=\iota_m'^*=\iota_m.
\end{split}
\end{equation*}
Then 
\[
\operatorname{Irr}(\mathscr{H}_{\bar{K}_1,B_n}(q)^{\dagger})=
\{\pi_{\lambda_i,\dagger},\pi_{\mu_j,\sharp,\dagger}^{\pm},
\pi_{\nu_k,\flat,\dagger}^{\pm},
\pi_{\kappa_l,\natural,\dagger}^{\pm},
\pi_{\iota_m,\dagger}^{++},\pi_{\iota_m,\dagger}^{+-},\,
\pi_{\iota_m,\dagger}^{-+},\,\pi_{\iota_m,\dagger}^{--}\}_{i,j,k,l,m}
\]
is a basic set of irreducible representations of $\mathscr{H}_{\bar{K}_1,B_n}(q)^{\dagger}$. 
Moreover $\mathscr{H}_{\bar{K}_1,B_n}(q)^{\dagger}$ is semisimple. 
\end{Theo}
\begin{proof}
We show that the elements of 
$\operatorname{Irr}(\mathscr{H}_{\bar{K}_1,B_n}(q)^{\dagger})$ are mutually inequivalent. 
We have already proved (1)--(4) below. \\
{\indent}(1) $\pi_{\mu_j,\sharp,\dagger}^+{\ncong}\pi_{\mu_j,\sharp,\dagger}^-$, {\qquad}
{\indent}(2) $\pi_{\nu_k,\flat,\dagger}^+{\ncong}\pi_{\nu_k,\flat,\dagger}^-$, {\qquad}
{\indent}(3) $\pi_{\kappa_l,\natural,\dagger}^+{\ncong}\pi_{\kappa_l,\natural,\dagger}^-$, \\
{\indent}(4) $\pi_{\iota_m,\dagger}^{++}$, $\pi_{\iota_m,\dagger}^{+-}$, 
$\pi_{\iota_m,\dagger}^{-+}$, $\pi_{\iota_m,\dagger}^{--}$ are mutually inequivalent. \\
Indeed, (1)--(3) have proved in Proposition 7.3, and (4) in Proposition 7.7. 
For brevity, we denote $\mathscr{H}_{\Bar{K}_1,B_n}(1,q)^{\sharp}$ by $\mathscr{H}^{\sharp}$ 
and $\mathscr{H}_{\Bar{K}_1,B_n}(1,q)^{\dagger}$ by $\mathscr{H}^{\dagger}$, and so on. 
We shall give inductions of elements of 
$\operatorname{Irr}(\mathscr{H}_{\bar{K}_1,B_n}(q)^{\dagger})$ to 
$\mathscr{H}=\mathscr{H}_{\Bar{K}_1,B_n}(1,q)$. We consider five cases depending upon 
the form of the $2$-tuple of Young diagrams. \\
case 1 : $\lambda_i{'}{\neq}\lambda_i,\lambda_i^*{\neq}\lambda_i,\lambda_i'^*{\neq}\lambda_i$\par
In this case, the induction of $V_{\lambda_i,\sharp,\dagger}$ to $\mathscr{H}^{\sharp}$ yields 
$(V_{\lambda_i,\sharp,\dagger})^{\mathscr{H}^{\sharp}}
{\cong}V_{\lambda_i,\sharp}{\oplus}V_{\lambda_i',\sharp}$ by Theorem 2.7 (1). 
Moreover, the induction of $V_{\lambda_i,\sharp}$ to $\mathscr{H}$ yields 
$(V_{\lambda_i,\sharp})^{\mathscr{H}}
{\cong}V_{\lambda_i}{\oplus}V_{\lambda_i^*}$ by Theorem 2.7 (1) again.
Thus we obtain 
$(V_{\lambda_i,\sharp,\dagger})^{\mathscr{H}}
{\cong}V_{\lambda_i}{\oplus}V_{\lambda_i^*}{\oplus}V_{\lambda_i'}{\oplus}V_{\lambda_i'^*}$. \\
case 2 : $\mu_j^*=\mu_j,\mu_j{'}{\neq}\mu_j,\mu_j'^*{\neq}\mu_j$\par
In this case, the induction of $V_{\mu_j,\sharp,\dagger}^+$ to $\mathscr{H}^{\sharp}$ yields 
$(V_{\mu_j,\sharp,\dagger}^+)^{\mathscr{H}^{\sharp}}
{\cong}V_{\mu_j,\sharp}^+{\oplus}V_{\mu_j',\sharp}^+$ by Proposition 7.7 (1) and Theorem 2.7 (1). 
Moreover, the induction of $V_{\mu_j,\sharp}^+$ to $\mathscr{H}$ yields 
$(V_{\mu_j,\sharp}^+)^{\mathscr{H}}{\cong}V_{\mu_j}$ by Theorem 2.7 (2). 
Thus we obtain $(V_{\mu_j,\sharp,\dagger}^+)^{\mathscr{H}}{\cong}V_{\mu_j}{\oplus}V_{\mu_j'}$. 
Similarly, $(V_{\mu_j,\sharp,\dagger}^-)^{\mathscr{H}}{\cong}V_{\mu_j}{\oplus}V_{\mu_j'}$ holds. \\
case 3 : $\nu_k^*{\neq}\nu_k,\nu_k{'}=\nu_k,\nu_k'^*{\neq}\nu_k$\par
Considering the induction of $V_{\nu_k,\sharp,\dagger}^+$ to $\mathscr{H}^{\flat}$, 
we obtain $(V_{\nu_k,\sharp,\dagger}^{\pm})^{\mathscr{H}}{\cong}V_{\nu_k}{\oplus}V_{\nu_k^*}$ 
in the same fashion as in the case 2. \\
case 4 : $\kappa_l^*{\neq}\kappa_l,\kappa_l{'}{\neq}\kappa_l,\kappa_l'^*=\kappa_l$\par
Considering the induction of $V_{\kappa_l,\sharp,\dagger}^+$ to $\mathscr{H}^{\natural}$, 
we obtain $(V_{\kappa_l,\sharp,\dagger}^{\pm})^{\mathscr{H}}{\cong}V_{\kappa_l}{\oplus}V_{\kappa_l'}$ 
in the same fashion as in the case 2. \\
case 5 : $\iota_m{'}=\iota_m^*=\iota_m'^*=\iota_m$\par
In this case, the induction of $V_{\iota_m,\sharp,\dagger}^{++}$ to $\mathscr{H}^{\sharp}$ yields 
$(V_{\iota_m,\sharp,\dagger}^{++})^{\mathscr{H}^{\sharp}}
{\cong}V_{\iota_m,\sharp}^+$ by Corollary 7.8 and Theorem 2.7 (2). 
Moreover, the induction of $V_{\iota_m,\sharp}^+$ to $\mathscr{H}$ yields 
$(V_{\iota_m,\sharp}^+)^{\mathscr{H}}{\cong}V_{\iota_m}$ by Theorem 2.7 (2). 
Thus we obtain $(V_{\iota_m,\sharp,\dagger}^{++})^{\mathscr{H}}{\cong}V_{\iota_m}$. 
Similarly, we get 
$(V_{\iota_m,\sharp,\dagger}^{+-})^{\mathscr{H}}{\cong}
(V_{\iota_m,\sharp,\dagger}^{-+})^{\mathscr{H}}{\cong}
(V_{\iota_m,\sharp,\dagger}^{--})^{\mathscr{H}}{\cong}V_{\iota_m}$. \par
Those induced representations are mutually non-isomorphic, hence we can conclude that 
the elements of $\operatorname{Irr}(\mathscr{H}_{\bar{K}_1,B_n}(q)^{\dagger})$ are 
mutually inequivalent. \par
Finally, we shall show the semisimplicity of $\mathscr{H}_{\Bar{K}_1,B_n}(1,q)^{\dagger}$. 
We define a map $\pi_n$ to be 
\begin{equation*}
\pi_n : x{\in}\mathscr{H}^{\dagger}{\longmapsto}
\bigoplus_{\pi{\in}\operatorname{Irr}(\mathscr{H}^{\dagger})}\pi(x)
\end{equation*} 
Then, by theorems of Burnside and Frobenius-Schur, 
$\mathscr{H}^{\dagger}$ has a quotient 
$\pi_n(\mathscr{H}^{\dagger})$ 
isomorphic to the semisimple algebra 
\[
\bigoplus_{\pi{\in}\operatorname{Irr}(\mathscr{H}^{\dagger})}
\operatorname{End}_{\Bar{K}_1}\Bar{K}_1^{{\deg}\pi}{\cong}
\bigoplus_{\pi{\in}\operatorname{Irr}(\mathscr{H}^{\dagger})}
\operatorname{Mat}({\deg}\pi,\Bar{K}_1). 
\]
We claim that this semisimple algebra has dimension $2^{n-2}n!$. 
On the other hand, $\dim_{\Bar{K}_1}\,\mathscr{H}^{\dagger}=2^{n-2}n!$, 
Therefore we can conclude that 
$\operatorname{Irr}(\mathscr{H}_{\bar{K}_1,B_n}(1,q)^{\dagger})$ 
is a basic set of irreducible representations and that 
$\mathscr{H}_{\bar{K}_1,B_n}(1,q)^{\dagger}$ is semisimple. 
\end{proof}

\end{document}